\documentclass[12pt]{amsart}
\usepackage{amssymb}
\usepackage{amsthm}
\usepackage{amscd}

\usepackage{amsmath}
\usepackage{epic,eepic}
\usepackage{graphicx}
\DeclareMathAlphabet{\mathpzc}{OT1}{pzc}{m}{it}

\usepackage{mathrsfs}
\usepackage{latexsym}

\usepackage{pifont}

\theoremstyle{plain}
\newtheorem{lemma}{Lemma}[section]
\newtheorem{prop}[lemma]{Proposition}
\newtheorem{thm}[lemma]{Theorem}
\newtheorem{cor}[lemma]{Corollary}
\newtheorem{aplemma}{Lemma~A.\hspace{-1.5mm}}
\newtheorem{approp}{Proposition~A.\hspace{-1.5mm}}
\newtheorem{apthm}{Theorem~A.\hspace{-1.5mm}}
\newtheorem{apcor}{Corollary~A.\hspace{-1.5mm}}
\newtheorem{intthm}{Theorem}

\usepackage[all]{xy}
\usepackage {cancel}

\newtheorem{conj}[lemma]{Conjecture}

\newcommand{\SSP}{\vspace{3mm}}
\newcommand{\LSP}{\vspace{5mm}}

\theoremstyle{definition}

\newtheorem{rema}[lemma]{Remark}

\newtheorem{remb}{Remark}

\newtheorem{defi}[lemma]{Definition}
\newtheorem{exa}[lemma]{Example}
\newtheorem{aprem}{Remark~A.\hspace{-1.5mm}}
\newtheorem{apdefi}{Definition~A.\hspace{-1.5mm}}
\newcommand{\bde}{\begin{defi}}
\newcommand{\ede}{\end{defi}\vspace{1mm}}
\newcommand{\ble}{\begin{lemma}}
\newcommand{\ele}{\end{lemma}}
\newcommand{\bpr}{\begin{prop}}
\newcommand{\epr}{\end{prop}}
\newcommand{\bt}{\begin{thm}}
\newcommand{\et}{\end{thm}}
\newcommand{\bco}{\begin{cor}}
\newcommand{\eco}{\end{cor}}
\newcommand{\bre}{\begin{rema}}
\newcommand{\ere}{\end{rema}}
\newcommand{\brea}{\begin{rema}}
\newcommand{\erea}{\end{rema}\vspace{1mm}}
\newcommand{\breb}{\begin{remb}}
\newcommand{\ereb}{\end{remb}\vspace{1mm}}
\newcommand{\bex}{\begin{exa}}
\newcommand{\eex}{\end{exa}}
\newcommand{\bpf}{\begin{proof}}
\newcommand{\epf}{\end{proof}\vspace{1mm}}

\newcommand{\bade}{\begin{apdefi}}
\newcommand{\eade}{\end{apdefi}}
\newcommand{\bale}{\begin{aplemma}}
\newcommand{\eale}{\end{aplemma}}
\newcommand{\bapr}{\begin{approp}}
\newcommand{\eapr}{\end{approp}}
\newcommand{\bat}{\begin{apthm}}
\newcommand{\eat}{\end{apthm}}
\newcommand{\baco}{\begin{apcor}}
\newcommand{\eaco}{\end{apcor}}
\newcommand{\bare}{\begin{aprem}}
\newcommand{\eare}{\end{aprem}}

\newcommand{\be}{\begin{enumerate}}
\newcommand{\ee}{\end{enumerate}}
\newcommand{\bcd}{\[\begin{CD}}
\newcommand{\ecd}{\end{CD}\]}
\newcommand{\bit}{\begin{itemize}}
\newcommand{\eit}{\end{itemize}}
\newcommand{\bq}{\begin{quote}}
\newcommand{\eq}{\end{quote}}
\newcommand{\ba}{\begin{array}}
\newcommand{\ea}{\end{array}}

\newcommand{\mcD}{\mathcal{D}}
\newcommand{\mcE}{\mathcal{E}}
\newcommand{\mcF}{\mathcal{F}}
\newcommand{\mcG}{\mathcal{G}}

\newcommand{\mcI}{\mathcal{I}}

\newcommand{\mcK}{\mathcal{K}}
\newcommand{\mcL}{\mathcal{L}}
\newcommand{\mcM}{\mathcal{M}}
\newcommand{\mcN}{\mathcal{N}}
\newcommand{\mcO}{\mathcal{O}}
\newcommand{\mcP}{\mathcal{P}}
\newcommand{\mcQ}{\mathcal{Q}}

\newcommand{\mcS}{\mathcal{S}}
\newcommand{\mcT}{\mathcal{T}}

\newcommand{\mcV}{\mathcal{V}}

\newcommand{\mbC}{\mathbb{C}}

\newcommand{\mbF}{\mathbb{F}}

\newcommand{\mbP}{\mathbb{P}}

\newcommand{\mbZ}{\mathbb{Z}}

\newcommand{\mfS}{\mathfrak{S}}

\newcommand{\mpf}{\mathpzc{f}}

\newcommand{\msF}{\mathscr{F}}

\newcommand{\vin}{\rotatebox{90}{$\subseteq$}}
\newcommand{\migi}{\rightarrow}
\newcommand{\longmigi}{\longrightarrow}

\newcommand{\isom}{\stackrel{\sim}{\migi}}

\newcommand{\migiincl}{\hookrightarrow}

\newcommand{\migisurj}{\twoheadrightarrow}

\newcommand{\mr}{\mathrm}
\newcommand{\hidden}[1]{\,}

\newcommand{\N}{N}

\newcommand{\DD}{L}
\newcommand{\OR}{m}
\newcommand{\DIM}{l}
\newcommand{\DEG}{d}

\newcommand{\Y}{X}

\begin{document}

\title[Dormant opers and Gauss maps in positive  characteristic]{Dormant opers and Gauss maps  \\ in positive  characteristic}
\author{Yasuhiro Wakabayashi}
\date{}
\markboth{Yasuhiro Wakabayashi}{}
\maketitle
\footnotetext{Y. Wakabayashi: Department of Mathematics, Tokyo Institute of Technology, 2-12-1 Ookayama, Meguro-ku, Tokyo 152-8551, JAPAN;}
\footnotetext{e-mail: {\tt wkbysh@math.titech.ac.jp};}
\begin{abstract}

The Gauss map of a given projective variety is the rational map that
sends a smooth point to the tangent space at that point, considered as a point of the Grassmann variety. The present paper aims to generalize a result by H. Kaji on Gauss maps in positive characteristic and establish an interaction with the study of dormant opers, as well as Frobenius-projective structures. We first prove a correspondence between dormant opers on a smooth projective variety $X$ and closed immersions from $X$ into a projective space with purely inseparable Gauss map. By using this, we determine the subfields of the function field of a smooth curve in positive characteristic induced by Gauss maps. Moreover, the correspondence gives us a Frobenius-projective structure on a Fermat hypersurface. This example embodies an exotic phenomenon of algebraic geometry in positive characteristic.

\end{abstract}
\tableofcontents 

\section*{Introduction}
\LSP

\subsection*{0.1} \label{S001}

Let $X$ be 
an algebraic  variety of dimension $\DIM >0$ over an algebraically closed field $k$  embedded in the  projective space $\mbP^\DD$ for some $\DD >0$.
Denote by $\mr{Grass} (\DIM +1, \DD +1)$ the Grassmann variety classifying $(\DIM +1)$-dimensional quotient spaces of the $k$-vector space $k^{\DD +1}$;
  it may be identified with the space of $\DIM$-planes in  $\mbP^\DD$.
The {\it Gauss map} is the rational morphism
$\gamma : X \dashrightarrow \mr{Grass} (\DIM +1, \DD +1)$ that assigns to a smooth point $x$ 
the embedded  tangent space to $X$ at $x$ in $\mbP^\DD$.
 
 The notion of Gauss map is generalized  (cf. \S\,\ref{SS01})
 by using linear spaces tangent to higher order, often called the osculating spaces; see, e.g.,  
    ~\cite{De}, ~\cite{Fra},    
    and ~\cite{Poh}.
 Also, a  different generalization of Gauss map can be found  in, e.g., ~\cite{Zak}.
 The study of the Gauss map and such generalizations  has been a subject of algebraic geometry for a long time.

  It is a well-known fact   that
 the Gauss map of a smooth non-linear subvariety of a projective space in characteristic $0$
  is 
 finite and birational onto its image (cf. ~\cite[(I, 2.8)]{Zak}).
On the other hand, when the base field has positive characteristic,
the birationality is no longer true in general,  and the Gauss map can be inseparable.
Various  properties  of Gauss maps in positive characteristic have been investigated by many mathematician;  see, e.g., ~\cite{FuKa1}, ~\cite{FuKa2}, ~\cite{Kaj1}, ~\cite{Kaj2}, ~\cite{KaNo}, ~\cite{KlePi},    ~\cite{Nom1}, and ~\cite{Nom2}.

For example, a result by H. Kaji
(cf. ~\cite[Corollary 6.2]{Kaj2})
 asserts that giving a closed immersion $X \migiincl \mbP^\DD$ from  a given smooth projective curve $X$  with  purely inseparable Gauss map of degree $p^N$ ($N >0$)
is equivalent to giving 
 a  rank $2$ vector bundle
  on the $N$-th Frobenius twist $X^{(N)}$  of $X$ satisfying certain conditions.
   As mentioned in Remark \ref{Rf99} of the present paper, this data 
 may be interpreted as  
  a {\it dormant $\mr{GL}_2$-oper of level $N$}, in the sense of ~\cite[Definition 4.2.1]{Wak9}; that is to say, 
  it gives  a certain rank $2$ vector  bundle on $X$ equipped with both an action of the sheaf  $\mcD_X^{(N-1)}$ (:= the ring of differential operators on $X$ of level $N-1$,  introduced in ~\cite{PBer1})
 and  a Hodge subbundle
 satisfying a strict form of Griffiths transversality.
 We refer the reader to, e.g.,  ~\cite{Mzk2},  ~\cite{Wak}
  for the study of dormant $\mr{GL}_2$-opers on curves (which are also known as {\it dormant indigenous bundles}), and the higher-rank cases were investigated  in, e.g.,  ~\cite{JP}, ~\cite{JRXY}, and ~\cite{Wak8}.
  
\LSP
\subsection*{0.2} \label{S003}

The present paper aims to refine and generalize  Kaji's result in order   to build an interaction   between the studies of dormant opers and  Gauss maps in positive characteristic.
To do this in a unified formulation involving multi-dimensional varieties, we 
introduce the notion of a  {\it dormant $(n, N)$-oper} (cf. Definition \ref{D45}),  which 
extends  the classical notion of a higher-level dormant  oper.
(However,  our discussion  deals essentially only with the case where 
$n=2$ or the underlying variety has dimension $1$.)

Let $X$ be a smooth projective variety   over an algebraically closed field $k$  of characteristic $p>2$ and $\chi := (n, N, d)$  a triple of positive  integers with $1 < n \leq p$.
We shall write
\begin{align}
\mr{Op}_{\chi, +\mr{imm}}^{^\mr{Zzz...}}
\end{align}
(cf. (\ref{E4509})) for the set of isomorphism classes of dormant $(n, N)$-opers  on $X$ equipped with certain additional data.
On the other hand,
we write 
\begin{align}
\mr{G au}_{\lambda}^{\mr{F}}
\end{align}
(cf. (\ref{e4872})), where $\lambda := (n-1, N, \DEG)$,
for the set of isomorphism classes of closed immersions $\iota : X \migiincl \mbP^\DD$ (for some $\DD >0$) of degree $d$ whose Gauss maps of order $n-1$ factor through the $N$-th relative Frobenius morphism.
Then, the main result in the first half of the present paper is the following assertion.

\SSP
\begin{intthm} [cf. Theorem \ref{T13} for the full statement] \label{TA}
Suppose that the quadruple  $(X, n,N, d)$ satisfies one of the conditions (a) and (b) described in \S\,\ref{SS71}.
Then, there exists a canonical injection of sets
\begin{align}
 \Xi_{\chi} : \mr{G au}_{\lambda}^{\mr{F}} \migiincl \mr{Op}_{\chi, + \mr{imm}}^{^\mr{Zzz...}}.
 \end{align}
 Moreover, this map is bijective when $n =2$.
 \end{intthm}
\SSP

\LSP
\subsection*{0.3} \label{S004}

We here describe two applications of   the above theorem proved in the second half of the present paper.
As the first application, we use the bijection $\Xi_{(2, \N, d)}$ to  determine the subfields of the function field of a given  curve induced by  Gauss maps.

Now, let $k$ be as above
and $X$   a smooth projective curve over $k$ of genus $g>1$.
Denote by 
$K(X)$ the function field of $X$ and by 
$\mcK$ the set of subfields $K$ of $K(X)$ satisfying the following condition:  There exists a closed immersion $\iota : X \migiincl \mbP^\DD$ (for some $\DD > 1$)
  such that 
the extension of function fields defined by the Gauss map associated to $(X, \iota)$ coincides with $K(X)/K$.
 
 H. Kaji proved  (cf. ~\cite[Corollaries 2.3 and  4.4]{Kaj2})
 that $K(X)$ itself belongs to $\mcK$, and that
  any subfield in  $\mcK$ is of the form $K(X)^{p^N} := \{ v^{p^N} \, | \, v \in K (X) \}$ for some integer $N\geq 0$, i.e.,  the inclusion relation $\mcK \subseteq \{ K (X)^{p^N} \, | \, N \geq 0 \}$ holds.

 To improve this result, we combine the above theorem with 
 the previous study of higher-level dormant  opers on curves developed in ~\cite{Wak6},
  in which we have shown the existence of a dormant $(2, N)$-opers  for every $N>0$.
The resulting assertion is described as follows.

\SSP
\begin{intthm} [= Theorem \ref{C091}] \label{TB}
Let us keep the above notation, and suppose that $2 < p$ and  $p \nmid (g-1)$.
Then, the following equality of sets holds:
\begin{align}
\mcK = \left\{ K(X)^{p^N} \, \Big| \, N \geq 0 \right\}.
\end{align}
 In particular, for every nonnegative integer $N$,
 there exists a closed immersion $X \migiincl \mbP^\DD$ (for some positive integer $\DD$) such that 
 the extension of function fields
 defined by the  Gauss map coincides with 
  $K(X)/K(X)^{p^N}$.
 \end{intthm}

\LSP
\subsection*{0.4} \label{S005}

The second application concerns higher-dimensional varieties.
We shall recall  from ~\cite[Definition 1.2.1]{Wak6} (or ~\cite[Definition 2.1]{Hos2}) the notion of 
an $F^N$-projective structure;
 this 
is a positive characteristic analogue of  the classical notion of a projective structure on a complex manifold discussed in, e.g., ~\cite{Gun}, ~\cite{KO1}, and ~\cite{KO2}.
Roughly speaking, an $F^N$-projective structure on a smooth variety $X$ (for $N>0$) is   a maximal collection of \'{e}tale coordinate charts on $X$ valued in $\mbP^{\mr{dim} (X)}$ 
 whose transition functions descend to the $N$-th Frobenius twist $X^{(N)}$ of $X$.

One ultimate goal of the study of $F^N$-projective structures is to give a  complete answer to  (the positive characteristic version of) the classification problem, starting with S. Kobayashi and T. Ochiai (cf. ~\cite{KO1}, ~\cite{KO2}), of varieties admitting projective structures.
In ~\cite{Wak6}, we  developed  the classification for some classes of varieties, including curves, surfaces, and Abelian varieties.
The difficulty is that, unlike the $1$-dimensional case,
  there are nontrivial obstructions for the existence of an $F^N$-projective structure. 
Indeed, because  of our lack of technical knowledge, only a few examples  have been previously found for higher dimensions.

However,  the bijection $\Xi_{(2, N, d)}$ asserted  in Theorem \ref{TA} enables us to  construct   an $F^N$-projective structure   by using an example  of a projective variety whose Gauss map is in a certain special situation. 
The assertion obtained in the present paper  is described  as follows.

\SSP
\begin{intthm} [= Theorem \ref{C01}] \label{TC}
Let $N$ be a positive integer and $\DD$ an integer with $\DD \neq 3$.
Denote by  $X$ the Fermat hypersurface of degree $p^N +1$ in $\mbP^\DD$ (cf. (\ref{E481})).
Then, $X$ admits an  $F^N$-projective structure 
\begin{align}
\mcS^\blacklozenge_{\mr{Gau}}
\end{align}
 arising from the Gauss map associated to the natural closed immersion $X \migiincl \mbP^\DD$.
Moreover, if $p\nmid \DD (\DD +1)$, then $X$ admits no $F^{2N +1}$-projective structures.
 \end{intthm}
\SSP

As mentioned in Remark \ref{R05}, the existence of such an $F^N$-projective structure may be thought of as an exotic phenomenon of algebraic geometry in positive characteristic. 
In fact, any unirational projective complex manifold  which is not isomorphic to projective spaces, such as Fermat hypersurfaces, admits no projective structures.

Also,
note that the only previous examples of $F^N$-projective structures on higher-dimensional varieties except for those on projective spaces 
were derived  from $F^N$-affine  structures on Abelian varieties or smooth curves equipped with a Tango structure.
By calculating Chern classes on the Fermat hypersurface $X$, we see that
$\mcS^\blacklozenge_{\mr{Gau}}$
 cannot be constructed in that way (cf. Remark \ref{R04}).
This means that the $F^N$-projective  structure asserted in the above theorem 
is essentially a new example.

\LSP
\subsection*{Notation and Conventions} 

Throughout the present paper, we
fix a prime number $p$ and 
an algebraically closed field $k$ of  characteristic $p$.

By a {\it variety (over $k$)}, we mean 
 a connected integral  scheme of finite type over $k$.
Moreover,  by a {\it curve},
 we mean a variety over $k$ of dimension $1$.
Unless stated otherwise, we will always be working  over $k$; for example,  products  of varieties will be taken over $k$, i.e., $X_1 \times X_2 := X_1 \times_k X_2$.

Let $X$ be a variety over $k$.
We shall write  $\Omega_X$ (resp., $\mcT_X$) for  the sheaf of $1$-forms (resp., the sheaf of vector fields) on $X$ over  $k$.
If $\mcV$ is a vector bundle (i.e., a locally free coherent sheaf) on $X$, then we denote by $\mbP (\mcV)$ the projective bundle over $X$ associated to $\mcV$.

Next, let $N$ be a positive integer.
We shall denote by $X^{(N)}$ the {\it $N$-th Frobenius twist} of $X$, i.e., the base-change of $X$ by the $p^N$-th power map $k \migi k$.
The {\it $N$-th relative Frobenius morphism} is
denoted by $F^{(N)}_{X/k} : X \migi X^{(N)}$.
When $N =1$, we write $F_{X/k}$ instead of $F_{X/k}^{(1)}$.
Also, we set $X^{(0)} := X$ and  $F_{X/k}^{(0)} = \mr{id}_X$ for simplicity.

Recall from ~\cite[\S\,2.2]{PBer1}
the sheaf of differential operators  $\mcD_X^{(N-1)} := \mcD_{X/\mr{Spec}(k)}^{(N-1)}$   on $X$ of level $N-1$, where $\mr{Spec}(k)$ is equipped with the trivial $(N-1)$-PD structure.
If $\nabla$ is  
 a left $\mcD_X^{(0)}$-action on an $\mcO_X$-module $\mcF$ extending its $\mcO_X$-module structure, then we will use the same notation to denote   the corresponding connection $\mcF \migi \Omega_X \otimes \mcF$.
Also, for a $\mcD_{X}^{(N-1)}$-action $\nabla$ on  $\mcF$, 
we shall write $\nabla^{(0)}$ for the $\mcD_X^{(0)}$-action (or equivalently, the connection) on $\mcF$ induced by $\nabla$ via the natural morphism $\mcD_X^{(0)} \migi \mcD_X^{(N-1)}$.

For an $\mcO_{X^{(N)}}$-module $\mcG$, there exists a canonical left $\mcD_X^{(N-1)}$-action
\begin{align} \label{E445}
\nabla_{\mcG, \mr{can}}^{(N-1)} : \mcD_X^{(N-1)} \migi \mcE nd_k (F_{X/k}^{(N)*}(\mcG))
\end{align}
on the pull-back $F_{X/k}^{(N)*}(\mcG)$ with vanishing $p$-$(N-1)$-curvature (cf. ~\cite[Definition 3.1.1 and Corollary 3.2.4]{LeQu}).
Given    an $\mcO_X$-module $\mcF$ and a left $\mcD_X^{(N-1)}$-action $\nabla$
on $\mcF$ extending its $\mcO_X$-module structure,
we shall write $\mcF^\nabla$
 for the subsheaf of $\mcF$ on which $\mcD_X^{(N-1)+}$ acts as zero, where $\mcD_X^{(N-1)+}$ denotes the kernel of the canonical projection $\mcD_X^{(N-1)} \migisurj \mcO_X$.
Note that $\mcF^\nabla$ may be regarded as  an $\mcO_{X^{(N)}}$-module   via the underlying homeomorphism of $F_{X/k}^{(N)}$.
The trivial $\mcD_{X}^{(\N-1)}$-action  on $\mcO_X$ will be denoted by $\nabla_{\mr{triv}}^{(N-1)}$.

Finally, for each positive   integer $\DD$, 
we denote the $\DD$-dimensional  projective space over $k$ by
\begin{align} \label{e49926}
\mbP^{\DD}  := \mr{Proj}(k[t_0, \cdots, t_\DD]) \left(= \{ [t_0: t_1: \cdots : t_{\DD}] \, | \, (t_0, \cdots, t_{\DD}) 
 \neq
 (0, \cdots, 0)\} \right).
 \end{align}

\LSP
\subsection*{Acknowledgements} 
We are grateful for the many constructive conversations we had with {\it algebraic varieties in positive characteristic}, who live in the world of mathematics!
Our work was partially supported by Grant-in-Aid for Scientific Research (KAKENHI No. 21K13770).

\vspace{10mm}
\section{Gauss maps  in positive characteristic} \label{S01}
\SSP

In this section, 
we recall the higher-order Gauss map associated to a closed subvariety of a projective space.
After that,  we will observe that, under a certain assumption,  the Gauss map induces an action of the ring of differential operators  on a jet bundle (cf. (\ref{E69})).

\LSP
\subsection{}
 \label{SS01}

Let  $X$ be  a smooth projective  variety over $k$ of dimension $\DIM >0$.
Denote by $\mcI$ the ideal sheaf  defining the diagonal in $X \times X$.
Also, for each $i =1, 2$, we denote by $\mr{pr}_i$ the $i$-th projection $X \times X \migi X$.

Let  us fix a line bundle  $\mcL$  on $X$ and  a nonnegative  integer $\OR$.
The  sheaf 
\begin{align} \label{E23}
J_\OR (\mcL) := \mr{pr}_{1*}(\mr{pr}_2^* (\mcL) \otimes \mcO_{X \times X}/\mcI^{\OR+1})
\end{align}
 forms a vector bundle on $X$ of rank $\binom{\DIM+ \OR}{\OR}$,  and it is 
 called 
  the {\bf $\OR$-jet bundle of $\mcL$}.
This sheaf is equipped with an $(\OR+1)$-step  decreasing filtration 
\begin{align} \label{E70}
\{ J_\OR (\mcL)^j \}_{j=0}^{\OR+1}
\end{align}
 given by putting $J_\OR (\mcL)^0 := J_\OR (\mcL)$ and  $J_\OR (\mcL)^j := \mr{Ker}\left(J_\OR (\mcL) \migisurj J_{j-1}(\mcL)  \right)$ ($j=1, \cdots, \OR+1$).
For each $j=0, \cdots, \OR$,  we have an isomorphism of $\mcO_X$-modules
\begin{align} \label{E77}
S^{j}(\Omega_X) \otimes \mcL \isom J_\OR (\mcL)^j/J_\OR (\mcL)^{j+1},
\end{align}
where $S^{j}(\Omega_X)$ denotes the $j$-th symmetric product  of $\Omega_X$ over $\mcO_X$.

Note that  $\mr{pr}_{1*}(\mr{pr}_2^* (\mcL))$ is canonically isomorphic to 
the vector bundle $H^0 (X, \mcL) \otimes_k \mcO_X$.
By  applying the functor $\mr{pr}_{1*}(\mr{pr}_2^* (\mcL) \otimes (-))$ to the quotient $\mcO_{X \times X} \migisurj \mcO_{X \times X}/\mcI^{\OR+1}$,
we obtain an $\mcO_X$-linear morphism
\begin{align} \label{E49}
H^0 (X, \mcL) \otimes_k \mcO_X \migisurj J_\OR (\mcL).
\end{align}

Next, suppose that we are given a closed immersion $\iota : X \migiincl \mbP^\DD$ for some positive integer $\DD$.
This induces  the composite
\begin{align} \label{E50}
\alpha_{\iota}^\OR  : \mcO_X^{\oplus (\DD +1)} \left( =H^0 (\mbP^\DD, \mcO_{\mbP^\DD}(1)) \otimes_k \mcO_{X}\right) 
&\migi H^0 (X, \iota^*(\mcO_{\mbP^\DD}(1))) \otimes_k \mcO_X  \\ &\xrightarrow{(\ref{E49})} J_\OR (\iota^*(\mcO_{\mbP^\DD}(1))), \notag
\end{align}
where the first arrow is the morphism induced  from the natural morphism 
$\mcO_{\mbP^\DD}(1) \migi \iota_*(\iota^*(\mcO_{\mbP^\DD}(1)))$.
Let $\mr{Grass}(\binom{\DIM+\OR}{\OR}, \DD +1)$ denote
 the Grassmann variety classifying  $\binom{\DIM+\OR}{\OR}$-dimensional    quotient spaces of  the $k$-vector space $k^{\DD +1}$.
If $U^\OR_\iota$ denotes the open locus of $X$ where 
$\alpha_\iota^\OR$ is surjective,
then the restriction of $\alpha_\iota^\OR$ to $U_\iota^\OR$ determines a morphism
\begin{align} \label{Er453}
\gamma^\OR_\iota : U^\OR_\iota \migi 
 \mr{Grass}(\textstyle{\binom{\DIM+\OR}{\OR}}, \DD +1).
 \end{align}
 We  call it  the {\bf Gauss map of order $\OR$} associated to $(X, \iota)$. 
 When $\OR =1$, the morphism $\gamma_\iota^\OR$ coincides with the Gauss map in the classical  sense (cf. Introduction).

Given a triple of positive integers $\lambda := (\OR, N, \DEG)$, we shall denote by 
\begin{align} \label{e4872}
\mr{G au}_{\lambda}^\mr{F}
\end{align}
the set of isomorphism classes of
closed immersions $\iota : X \migiincl \mbP^{\DD}$ (for some $\DD >0$) satisfying the following two conditions:
\begin{itemize}
\item
The closed subvariety $\mr{Im}(\iota)$ of $\mbP^\DD$ has  degree $\DEG$;
\item
$U_\iota^m = X$ and $\gamma_\iota^\OR$ factors through $F_{X/k}^{(N)}$.
\end{itemize}
Here, two such closed immersions $\iota_i : X \migiincl \mbP^{\DD_i}$ ($i=1,2$) are said to be {\it isomorphic} if there exists an isomorphism 
$h : \mbP^{L_1} \isom \mbP^{L_2}$ 
  satisfying $\iota_2 = h\circ \iota_1$.
   Thus, it makes sense to speak of  the {\it isomorphism class} of a  closed immersion  $\iota : X \migiincl \mbP^{\DD}$ as above.

By putting  $\lambda_{N'} := (\OR, N', \DEG)$ for each positive integer $N'$, we obtain 
 the following sequence of  inclusions:
\begin{align} \label{ER556}
\mr{G au}_{\lambda_1}^\mr{F} \supseteq \mr{G au}_{\lambda_2}^\mr{F} \supseteq \mr{G au}_{\lambda_3}^\mr{F} \supseteq \cdots \supseteq \mr{G au}_{\lambda_N}^\mr{F} \supseteq  \cdots.
\end{align}

\LSP
\subsection{}
 \label{SS07}

Let    $\iota : X \migiincl  \mbP^{\DD}$ (where $\DD >0$) be a closed immersion,
 $\OR$ a nonnegative integer, and $N$  a positive integer.
  Suppose that
$U^\OR_\iota  = X$
and $\gamma_\iota^\OR$  factors through $F_{X/k}^{(N)}$.
Then, we can find a unique morphism
 $\breve{\gamma} : X^{(N)} \migi \mr{Grass}(\textstyle{\binom{\DIM+\OR}{\OR}}, \DD +1)$
with  $\breve{\gamma} \circ F_{X/k}^{(N)}  = \gamma_\iota^\OR$.
Let us  denote  the universal quotient on $\mr{Grass}(\textstyle{\binom{\DIM+\OR}{\OR}}, \DD +1)$  by 
\begin{align} \label{e37}
q_\mr{univ} : \mcO_{\mr{univ}}^{\oplus (\DD +1)} \migisurj \mcQ_{\mr{univ}},
\end{align}
where $\mcO_{\mr{univ}}$ denotes the structure sheaf.
The pull-back 
 of $q_{\mr{univ}}$ by $\breve{\gamma}$ defines an $\mcO_{X^{(N)}}$-linear surjection
$q_0 : 
\mcO_{X^{(N)}}^{\oplus (\DD +1)} \migisurj \mcQ_{0}$.
It follows from the definition of $\gamma_\iota^\OR$ that
there exists a unique isomorphism 
$\tau : F_{X/k}^{(N)*} (\mcQ_0) \isom J_\OR (\iota^* (\mcO_{\mbP^\DD}(1)))$ which makes the following diagram commute:
\begin{align} \label{E456}
\vcenter{\xymatrix@C=46pt@R=36pt{
&
 \mcO_{X}^{\oplus (\DD +1)}\ar[ld]_{F_{X/k}^{(N)*}(q_0)} \ar[rd]^{\alpha_\iota^\OR}&
\\
F_{X/k}^{(N)*} (\mcQ_0)\ar[rr]_{\tau}^{\sim} && J_\OR (\iota^* (\mcO_{\mbP^\DD}(1))).
}}
\end{align}
The $\mcD_X^{(N)}$-action $\nabla_{\mcQ_0, \mr{can}}^{(N-1)}$ (cf. (\ref{E445})) corresponds, via $\tau$, 
a
 $\mcD_{X}^{(N -1)}$-action
\begin{align} \label{E69}
\nabla_{\iota, \mr{Gau}}^{(N-1)} : \mcD_X^{(N-1)} \migi \mcE nd_{k} (J_\OR (\iota^* (\mcO_{\mbP^\DD}(1)))) 
\end{align}
on the $\OR$-th jet bundle $J_\OR (\iota^* (\mcO_{\mbP^\DD}(1)))$ of $\iota^* (\mcO_{\mbP^\DD}(1))$.
By definition,  the $\mcD_{X}^{(N -1)}$-action $\nabla_{\iota, \mr{Gau}}^{(N-1)}$ has  vanishing $p$-$(N-1)$-curvature.

\vspace{10mm}
\section{Dormant $(n, N)$-opers on a variety} \label{S06}
\SSP

In this section, the classical definition of a (higher-level) dormant oper is generalized  to multi-dimensional varieties.
The main result of this section describes a relationship between higher-order Gauss maps and  generalized dormant opers (cf. Theorem \ref{T13}).

\LSP
\subsection{}
\label{SS44}

Let $X$ be a smooth projective  variety  over $k$ of dimension $\DIM>0$, and let $N$, $n$ be  two integers with $N >0$, $p \geq  n>1$.
For a rank $n$ vector bundle $\mcV$ on $X$ and an integer $a$ with $1 \leq a < p$, we shall write
$T^a (\mcV)$ (cf. ~\cite[Definition 3.4]{Sun1}) for the subbundle of $\mcV^{\otimes a}$ 
generated locally by various sections $\sum_{\sigma \in \mfS_a} \breve{e}_{\sigma (1)} \otimes \cdots \breve{e}_{\sigma (a)}$, where $\mfS_a$ denotes the symmetric group of $a$ letters and  each $\breve{e}_i$ ($i=1, \cdots, a$) is an element in a fixed local  basis $\{ e_1, \cdots, e_n \}$ of $\mcV$.
Also, we set $T^0 (\mcV) := \mcO_X$.
Note that the subbundle $T^a (\mcV)$ does not depend on the choice of the local basis $\{e_1, \cdots, e_n\}$, and that since $a < p$  it is  isomorphic to the $a$-th symmetric product $S^a (\mcV)$ of $\mcV$ over $\mcO_X$ via the composite of natural morphisms $T^a (\mcV) \migiincl \mcV^{\otimes a} \migisurj S^a (\mcV)$.

Let us consider a collection of data
\begin{align}
\msF^\heartsuit := (\mcF, \nabla, \{ \mcF^j \}_{j=0}^n)
\end{align}
consisting of a vector bundle $\mcF$ on $X$, a $\mcD_X^{(N-1)}$-action $\nabla$ on $\mcF$ extending its $\mcO_X$-module structure, and an $n$-step decreasing filtration $\{ \mcF^j \}_{j=0}^n$ on $\mcF$ such that $\mcF^0 = \mcF$, $\mcF^n = 0$, and  $\mcF/\mcF^1$ is a line bundle.

\SSP
\bde \label{D45}
\begin{itemize}
\item[(i)]
We say that $\msF^\heartsuit$ is an {\bf  $(n, N)$-oper}
  on $X$ if it satisfies the following conditions:
\begin{itemize}
\item
 For each $j =1, \cdots, n-1$, the inclusion relation $\nabla^{(0)}(\mcF^j) \subseteq  \Omega_X\otimes \mcF^{j-1}$ holds and the {\it $\mcO_X$-linear} morphism 
 \begin{align} \label{e100}
\mr{KS}_{\msF^\heartsuit}^j :  \mcF^j/\mcF^{j+1} \migi \Omega_X \otimes (\mcF^{j-1}/\mcF^j)
 \end{align}
  induced naturally by $\nabla^{(0)}$  is injective.
  We call  $\mr{KS}_{\msF^\heartsuit}^j$ the {\it $j$-th Kodaira-Spencer map} associated to $\msF^\heartsuit$.
\item
For each $j=0, \cdots, n-1$, the $\mcO_X$-linear morphism
\begin{align} \label{E411}
\mcF^j/\mcF^{j+1} \migi \Omega_X^{\otimes j} \otimes (\mcF/\mcF^1)
\end{align}
obtained by composing various $\mr{KS}_{\msF^\heartsuit}^{j'}$'s has image $T^j (\Omega_X)  \otimes (\mcF/\mcF^1) \subseteq \Omega_X^{\otimes j}  \otimes (\mcF/\mcF^1)$.
(If  either $\DIM =1$ or $n=2$ is satisfied, then this condition is equivalent to the condition that (\ref{E411}) is an isomorphism for every $j$.)
\end{itemize} 
Moreover,  the notion of an isomorphism between two $(n, N)$-opers  can be defined in a natural manner, so we will omit the details of that definition.
  \item[(ii)]
  An $(n, N)$-oper $\msF^\heartsuit := (\mcF, \nabla, \{ \mcF^j \}_{j=0}^n)$  is called {\bf dormant} if $\nabla$ has vanishing $p$-$(N-1)$-curvature.
  \end{itemize}
\ede
\SSP

\begin{rema} \label{R90}
If the underlying variety $X$ has dimension $1$, then
the definition of an oper using the ring of   higher-level differential operators
    can be found in 
 ~\cite[Definition 4.2.1]{Wak9}.
In the cace  of  $n=2$ (but $X$ is arbitrary),
an equivalent definition  was discussed in  ~\cite[Definition 2.3.1]{Wak6} under the name of {\it (dormant) indigenous $\mcD_X^{(N-1)}$-module}.
Also, it follows from  ~\cite[Theorem A]{Wak6} that, when $p \nmid \mr{dim}(X)+1$,  
equivalence classes (with respect to a certain equivalence relation) of dormant $(2, N)$-opers  are in bijection with what we call
 {\it $F^N$-projective structures}.
 (In \S\,\ref{SS03}, we will mention briefly the definition of an $F^N$-projective structure.)
\end{rema} 
\SSP

\begin{rema} \label{R90f}
Here,  let us recall a typical example of a dormant $(p, 1)$-oper  on a multi-dimensional variety  provided by  X. Sun (cf.  ~\cite[Theorem 3.7]{Sun1}); this is  a generalization of the dormant $(p, 1)$-oper on a curve  discussed in ~\cite[\S\,5]{JRXY}.

Given a line bundle $\mcL$ on $X$, we shall set $\mcV := F_{X/k}^* (F_{X/k*}(\mcL))$.
The $\mcO_X$-module  $\mcV$ forms a vector bundle on $X$ and admits a connection $\nabla := \nabla^{(0)}_{F_{X/k*}(\mcL), \mr{can}}$ with vanishing $p$-curvature.
This sheaf  is equipped with a $p$-step decreasing filtration $\{ \mcV^j \}_{j=0}^p$ given by the following construction:
\begin{itemize}
\item
$\mcV^0 := \mcV$ and  $\mcV^1$ is the kernel of the morphism $\mcV \migisurj \mcL$ corresponding to the identity morphism $\mr{id}_{F_{X/k*}(\mcL)}$ via the adjunction relation ``$F_{X/k}^* (-) \dashv F_{X/k*}(-)$".
\item
For each $j=2, \cdots, p$, we define  $\mcV^{j}$ inductively as follows:
\begin{align}
\mcV^{j} := \mr{Ker} \left(\mcV^{j-1} \xrightarrow{\nabla} \Omega_X \otimes \mcV^{j-2} \xrightarrow{\mr{quotient}} \Omega_X \otimes (\mcV^{j-2}/\mcV^{j-1}) \right).
\end{align}
\end{itemize}
Then, 
 the resulting collection  $(\mcV, \nabla, \{ \mcV^j \}_{j=0}^p)$ forms a dormant $(p, 1)$-oper on $X$.
\end{rema} 
\SSP

\begin{rema} \label{Rf99}
Suppose that $X$ is a smooth projective curve of genus $g >1$ and $N$ is a positive integer.
According to  ~\cite[Corollary 6.2]{Kaj2},
 the existence of a closed immersion $X \migiincl \mbP^{\DD}$ with purely inseparable Gauss map of degree $p^N$  
is equivalent to the existence of a rank $2$ vector bundle $\mcG$ on $X^{(N)}$ satisfying the  following condition:
\begin{itemize}
\item[$(*)_\mcG$] : 
$F^{(N-1)*}_{X^{(1)}/k}(\mcG)$ is stable and $F^{(N)*}_{X/k}(\mcG) \cong J_1 (\mcL)$ for some line bundle $\mcL$ on $X$.
\end{itemize}
In this Remark,  we shall  examine the relationship between such vector bundles $\mcG$ and dormant $(2, N)$-opers.
Let us take   a dormant $(2, N)$-oper   $\msF^\heartsuit := (\mcF, \nabla, \{ \mcF^j \}_{j=0}^2)$ on $X$.

First, we shall prove the claim that, {\it for every  positive integer $N'$ with $N' \leq  N$,   the rank $2$ vector bundle $F_{X^{(N')}/k}^{(N-N')*}(\mcF^\nabla)$ on $X^{(N')}$ 
 is   stable}.
 (The following discussion is available for a general $n$, but we focus on the rank $2$ case for simplicity.)
Suppose, on the contrary, that $F_{X^{(N')}/k}^{(N-N')*}(\mcF^\nabla)$ is unstable, i.e.,  there exists a line subbundle $\mcL$ of $F_{X^{(N')}/k}^{(N-N')*}(\mcF^\nabla)$ of degree $\geq \frac{1}{2} \cdot a$, where $a := \mr{deg}(F_{X^{(N')}/k}^{(N-N')*}(\mcF^\nabla))$.
Then, $F^{(N')*}_{X/k}(\mcL) = p^{N'} \cdot \mr{deg}(\mcL) \geq \frac{p^{N'}}{2} \cdot a$.
Since $\mr{KS}^1_{\msF^\heartsuit}$ is an isomorphism, we have
\begin{align} \label{E0004}
\mr{deg}(\mcF^1) - \mr{deg}(\mcF/\mcF^1) = \mr{deg}(\Omega_X \otimes (\mcF/\mcF^1))-\mr{deg} (\mcF/\mcF^1) =2g-2.
\end{align}
On the other hand, the following equalities holds:
\begin{align} \label{E0005}
\mr{deg}(\mcF^1) + \mr{deg} (\mcF/ \mcF^1) = \mr{deg} (\mcF) = \mr{deg} (F_{X/k}^{(N')*}(F_{X^{(N')}/k}^{(N-N')*} (\mcF^\nabla))) = p^{N'} \cdot a,
\end{align}
where the second equality follows from the isomorphism 
\begin{align} \label{E4980}
\left(F_{X/k}^{(N')*}(F_{X^{(N')}/k}^{(N-N')*}(\mcF^\nabla)) =  \right) F^{(N)*}_{X/k}(\mcF^\nabla) \isom \mcF
\end{align}
 resulting from  ~\cite[Corollary 3.2.4]{LeQu}.
It follows from  (\ref{E0004}) and (\ref{E0005}) that $\mr{deg}(\mcF^1) = \frac{p^{N'}}{2} \cdot a + g-1$ and $\mr{deg}(\mcF/\mcF^1) = \frac{p^{N'}}{2} \cdot a - g+1$.
By comparing the respective degrees of $F_{X/k}^{(N')*}(\mcL)$ and $\mcF/\mcF^1$, we see that
 the composite of natural morphisms
 \begin{align} \label{E442}
 F^{(N')*}_{X/k}(\mcL) \migiincl  F_{X/k}^{(N')*}(F_{X^{(N')}/k}^{(N-N')*}(\mcF^\nabla))
   \xrightarrow{(\ref{E4980})} \mcF \migisurj \mcF/\mcF^1
 \end{align}
  coincides with the zero map.
  Hence, the composite of the first two morphisms in (\ref{E442}), which we shall denote by $h$,   factors through the inclusion $\mcF^1 \migiincl \mcF$.
  The resulting morphism $F_{X/k}^{(N')*}(\mcL) \migi \mcF^1$ is injective and an isomorphism at the generic point $\eta$ of $X$.
 Let us identify each local section of $F_{X/k}^{(N')*}(\mcL)$ with its image via this injection.
  Then,  the stalk of $\mcF^1$ at  $\eta$ is closed under the connection $\nabla_{\mcL, \mr{can}}^{(0)}$. 
 Since $\nabla_{\mcL, \mr{can}}^{(0)}$ is compatible with $\nabla$ via $h$,  the stalk of $\mcF^1$  at $\eta$ is also closed under $\nabla$.
  But,  it contradicts the assumption that $\mr{KS}_{\msF^\heartsuit}^1$ is an isomorphism.
  Consequently, $F_{X^{(N')}/k}^{(N-N')*}(\mcF^\nabla)$ turns out to be stable, and  the proof of the claim is completed.
  In particular, $F_{X^{(1)}/k}^{(N-1)*}(\mcF^\nabla)$ is stable.

Next, observe that,  since $\mr{KS}_{\msF^\heartsuit}^1$ is an isomorphism,
the following composite is an isomorphism:
\begin{align}
\mcD_{X, \leq 1}^{(N-1)}\otimes (\mcF/\mcF^1)^\vee \xrightarrow{\mr{inclusion}} \mcD_X^{(N-1)} \otimes \mcF^\vee \migi \mcF^\vee,
\end{align}
where $\mcD_{X, \leq 1}^{(N-1)}$ denotes the subsheaf of $\mcD_{X}^{(N-1)}$ consisting of differential operators of order $\leq 1$ and the second arrow denotes  the morphism induced naturally by the dual of $\nabla$.
If we write $\mcL := \mcF/\mcF^1$, then  the dual of this composite determines an isomorphism $\mcF \isom J_1 (\mcL)$. 
(This isomorphism preserves the filtration, i.e., restricts to an isomorphism $\mcF^1 \isom J_1 (\mcL)^1$.)
Thus, we conclude that the rank $2$ vector bundle $\mcF^\nabla$ on $X^{(N)}$ satisfies the condition $(*)_{\mcF^\nabla}$ described above.

One may verify that {\it the resulting assignment $\msF^\heartsuit \mapsto \mcF^\nabla$ gives a bijective correspondence between dormant $(2, N)$-opers on $X$ and rank $2$ vector bundles $\mcG$ on $X^{(N)}$ satisfying  $(*)_\mcG$.}
\end{rema} 

\LSP
\subsection{}
 \label{SS79}

We shall  prove  the following assertion concerning dormant $(n, N)$-opers.

\SSP
\bpr \label{L23}
Let $\msF^\heartsuit := (\mcF, \nabla, \{ \mcF^j \}_j)$ be a dormant $(n, N)$-oper  on $X$.
Denote by $\sigma$ the morphism $X \migi \mbP (\mcF)$ induced, via projectivization,  from the natural quotient $\mcF \migi \mcF/\mcF^1$.
Also, denote by 
 $\overline{\sigma}$  the  composite
\begin{align} \label{E48}
\overline{\sigma} : X \xrightarrow{\sigma}\mbP (\mcF) \left(= X \times_{X^{(N)}} \mbP (\mcF^\nabla) \right) \xrightarrow{\mr{projection}} \mbP (\mcF^\nabla).
\end{align}
Then, the following assertions hold:
\begin{itemize}
\item[(i)]
Suppose that $\mr{dim} (X) = 1$.
Then, $\overline{\sigma}$ is birational onto its image.
\item[(ii)] 
Suppose that $n=2$ (but $\mr{dim} (X)$ is arbitrary).
Then,  $\overline{\sigma}$ is
 a closed immersion.
\end{itemize}
\epr
\begin{proof}
First, we shall consider assertion (i).
Denote by $Y$ the nomalization of the image $\mr{Im}(\overline{\sigma})$ of $\overline{\sigma}$.
Suppose that the field extension $K (X)/K (Y)$ is nontrivial.
Since $F_{X/k}^{(N)} : X \migi X^{(N)}$ factors through a morphism $h : X \migi Y$,
there exists an integer $M$ with  $1 \leq M \leq N$ such that $Y = X^{(M)}$ and $h = F_{X/k}^{(M)}$.
Denote by $\mcG$ the pull-back of $\mcF^\nabla$ to $X^{(M)}$.
In particular, its pull-back $F^{(M)*}_{X/k}(\mcG)$ may be canonically identified with $\mcF$.
The section  $\left( Y = \right) X^{(M)} \migi \mbP (\mcG)$ induced by  the composite of 
the normalization $Y \migi \mr{Im}(\overline{\sigma})$ and the inclusion $\mr{Im}(\overline{\sigma}) \migiincl \mbP (\mcF^\nabla)$
determines a surjection $\mcG \migisurj \mcQ$ for some line bundle $\mcQ$ on $X^{(M)}$.
Under the identification $F^{(M)*}_{X/k}(\mcG) = \mcF$,
the pull-back of this surjection $\mcG \migisurj \mcQ$  by $F_{X/k}^{(M)}$ coincides  with the natural projection   $\mcF \migisurj \mcF/\mcF^1$, and the equality $\nabla^{(0)} = (\nabla_{\mcG, \mr{can}}^{(M-1)})^{(0)}$ holds.
Hence, the subbundle $\mcF^1$ of $\mcF$ (which may be identified  with the pull-back of $\mr{Ker}\left(\mcG \migisurj \mcQ \right)$) is closed  under $\nabla^{(0)}$.
But, this contradicts the fact that 
$\mr{KS}_{\msF^\heartsuit}^1$ is nonzero.
Thus, we have  $K (X) = K (Y)$, meaning  that
$\overline{\sigma}$
  is birational onto its image.
  This completes the proof of assertion (i).

Next, we shall prove assertion (ii).
For each point $x$ of $X$, we can find an open neighborhood  $U$ of $x$ in $X$ such that
there exists an $\mcO_{U^{(M)}}$-linear  isomorphism $\mcF^\nabla |_{U^{(N)}}\isom \mcO_{U^{(N)}}^{\oplus \DIM +1}$.
Let  us consider  the isomorphism  of $\mbP^\DIM$-bundles 
\begin{align} \label{e584}
\mbP (\mcF^\nabla |_{U^{(N)}}) \isom U^{(N)}\times \mbP^\DIM
\end{align}
 induced by this isomorphism. 
Since
$\mr{KS}_{\msF^\heartsuit}^1$ is an isomorphism,
the composite
\begin{align}
U \xrightarrow{\overline{\sigma}|_{U}} \mbP (\mcF^\nabla |_{U^{(N)}}) \xrightarrow{(\ref{e584})} U^{(N)} \times \mbP^\DIM \xrightarrow{\mr{projection}}
\mbP^\DIM
 \end{align}
is \'{e}tale (cf. ~\cite[Corollary 1.6.2]{Wak6}).
 This implies  that the restriction  $\overline{\sigma} |_U$ of $\overline{\sigma}$ is unramified.
By  applying this argument to  various points $x$ of $X$,   we see  that  $\overline{\sigma}$ is unramified.
Moreover, since $F_{X/k}^{(N)}$ factors through $\overline{\sigma}$,
 the morphism 
 $\overline{\sigma}$ is universally injective.
 Thus, $\overline{\sigma}$ turns out to be a  closed immersion.
 This completes the proof of  assertion (ii).
 \end{proof}
\SSP

Given a triple of  positive integers $\chi := (n, N, \DEG)$ with $1< n \leq  p$, we shall write 
\begin{align} \label{E4509}
\mr{Op}_{\chi, +\mr{bir}}^{^\mr{Zzz...}} \ \left(\text{resp.,} \ \mr{Op}_{\chi, +\mr{imm}}^{^\mr{Zzz...}}\right)
\end{align}
for the set of isomorphism classes of 
 pairs $\mpf := (\msF^\heartsuit, q)$ consisting of  a dormant $(n, N)$-oper $\msF^\heartsuit := (\mcF, \nabla, \{ \mcF^j \}_j)$ on $X$ 
  and
 a  surjective morphism of $\mcD_{X}^{(\N-1)}$-modules $q : (\mcO_X, \nabla_\mr{triv}^{(N-1)})^{\oplus (\DD_\mpf +1)}$ $\migisurj (\mcF, \nabla)$ for some $\DD_\mpf >0$
  that satisfies the following two conditions:
 \begin{itemize}
  \item
 The morphism 
 \begin{align} \label{E4444}
 \iota_{\mpf} : X \migi \mbP^{\DD_\mpf}
 \end{align}
  determined, via projectivization,  by the composite of $q$ and the surjection $\mcF \migisurj \mcF/\mcF^1$ is
 birational onto its image (resp., a closed immersion);
 \item
The degree of the closed subvariety $\mr{Im}(\iota_{\mpf})$ of $\mbP^{\DD_\mpf}$
  is equal to  $\DEG$.
  (This is equivalent to the condition that 
    $\mcF/\mcF^1$ has degree $\DEG$ with respect to the ample line bundle   $\iota_{\mpf}^*(\mcO_{\mbP^{\DD_\mpf}}(1))$, in the sense of ~\cite[Definition 1.2.11]{HL}.)
 \end{itemize}
Here,  two such pairs $\mpf_i:= (\msF^\heartsuit_i, q_i)$ ($i =1,2$)
 are said to be {\it isomorphic} if there exists a pair $(h_\msF, h_\mcO)$ consisting of an isomorphism of $(n, N)$-opers $h_\msF : \msF_1^\heartsuit \isom \msF_2^\heartsuit$ and an isomorphism  of  $\mcD_X^{(\N -1)}$-modules $h_\mcO :  (\mcO_X, \nabla_\mr{triv}^{(N-1)})^{\oplus (\DD_{\mpf_1} +1)} \isom (\mcO_X, \nabla_\mr{triv}^{(N-1)})^{\oplus (\DD_{\mpf_2}+1)}$ satisfying $q_2 \circ h_\mcO = h_\msF \circ q_1$.
Thus, we can define the {\it isomorphism class} of a pair $\mpf := (\msF^\heartsuit, q)$ as above.
It is clear that $\mr{Op}_{\chi, +\mr{imm}}^{^\mr{Zzz...}} \subseteq \mr{Op}_{\chi, +\mr{bir}}^{^\mr{Zzz...}}$.

\LSP
\subsection{}
 \label{SS79ff}
 Hereinafter, we shall use the notation $\Box$ to  denote either  ``$\mr{bir}$" or ``$\mr{imm}$".
Let $\mpf := (\msF^\heartsuit, q)$ (where $\msF^\heartsuit := (\mcF, \nabla, \{ \mcF^j \}_j)$)
  be a pair classified by $\mr{Op}_{\chi, + \Box}^{^\mr{Zzz...}}$.
For a positive integer $N'$ with $N'<N$, 
denote by $\nabla^{(\N'-1)}$ the $\mcD_X^{(\N' -1)}$-action on $\mcF$ induced from $\nabla$ via the natural morphism $\mcD_X^{(N' -1)} \migi \mcD_X^{(N-1)}$.
 Then, the collection  $\msF^{\heartsuit (\N')} := (\mcF, \nabla^{(\N'-1)}, \{ \mcF^j \}_j)$ forms a dormant $(n, N')$-oper and the morphism $q$ determines  
 a surjective morphism of $\mcD_{X}^{(N' -1)}$-modules $(\mcO_X, \nabla_\mr{triv}^{(N'-1)})^{\oplus (\DD_\mpf +1)} \migisurj (\mcF, \nabla^{(N'-1)})$.
 In particular, the pair $\mpf^{(N')} := (\msF^{\heartsuit (N')}, q)$ is an element of $\mr{Op}_{\chi, + \Box}^{^\mr{Zzz...}}$, where $\chi' := (n, N', \DEG)$.
 The map of sets $\mr{Op}_{\chi, + \Box}^{^\mr{Zzz...}}\migi \mr{Op}_{\chi', + \Box}^{^\mr{Zzz...}}$ given by 
  $\mpf \mapsto \mpf^{(N')}$  is  verified to be injective.
 This injection allows us to  regard  $\mr{Op}_{\chi, + \Box}^{^\mr{Zzz...}}$ as a subset of $\mr{Op}_{\chi', + \Box}^{^\mr{Zzz...}}$.

Thus, 
by putting  $\chi_{N'} := (n, N', \DEG)$ for each positive integer $N'$, we obtain 
 the following diagram  of  inclusions:
\begin{align} \label{ER557}
&\mr{Op}_{\chi_1, +\mr{bir}}^{^\mr{Zzz...}}
  \  \supseteq \ \mr{Op}_{\chi_2, +\mr{bir}}^{^\mr{Zzz...}}
\ \supseteq \mr{Op}_{\chi_3, +\mr{bir}}^{^\mr{Zzz...}}
  \ \,\supseteq \cdots \, \supseteq \,
  \mr{Op}_{\chi_N, +\mr{bir}}^{^\mr{Zzz...}}  \   \supseteq  \cdots \\
  & \hspace{5mm}\vin \hspace{20mm} \vin \hspace{20mm} \vin  \hspace{15mm} \cdots \hspace{10mm} \vin \hspace{17mm}  \cdots \notag \\
&\mr{Op}_{\chi_1, +\mr{imm}}^{^\mr{Zzz...}}
 \supseteq \mr{Op}_{\chi_2, +\mr{imm}}^{^\mr{Zzz...}}
 \supseteq \mr{Op}_{\chi_3, +\mr{imm}}^{^\mr{Zzz...}}
  \supseteq \cdots \, \supseteq 
  \,\mr{Op}_{\chi_N, +\mr{imm}}^{^\mr{Zzz...}}
  \supseteq  \cdots. \notag
\end{align}

The following assertion gives a necessary condition for the set $\mr{Op}_{\chi, +\mr{bir}}^{^\mr{Zzz...}}$ being nonempty.
\SSP
\bpr \label{Dfe2}
Suppose that there exists a pair  $\mpf := (\msF^\heartsuit, q)$ 
classified by
 $\mr{Op}_{\chi, +\mr{bir}}^{^\mr{Zzz...}}$.
 Denote by 
 $\mr{deg}(\Omega_X)$ the degree of $\Omega_X$ with respect to 
 the ample divisor $H$ determined by $\iota_\mpf^*(\mcO_{\mbP^{\DD_\mpf}}(1))$, i.e., 
   $\mr{deg}(\Omega_X) := c_1 (\Omega_X) \cdot H^{\DIM -1}$.
Then, we have
\begin{align} \label{Eeer}
\frac{1}{p^N} \cdot  \left(\mr{deg}(\Omega_X) + \frac{\DEG (\DIM +1)}{n-1} \right) \cdot \binom{\DIM + n-1}{n-2} \in \mbZ.
\end{align}
In particular, if 
$X$ is a   smooth projective curve of genus $g$, 
then  $\mr{Op}_{\chi, +\mr{bir}}^{^\mr{Zzz...}} = \mr{Op}_{\chi, +\mr{imm}}^{^\mr{Zzz...}}  = \emptyset$ unless  $p^N \mid n ((n-1)(g-1)+ \DEG)$.
\epr
\begin{proof}
Let   $(\mcF, \nabla, \{ \mcF^j \}_j)$ be the collection of data defining $\msF^\heartsuit$.
It follows from ~\cite[Lemma 4.3]{Sun1} that
for each integer $j$ with $0 \leq j \leq n-1 \left(< p \right)$, the following equalities hold:
\begin{align} \label{FG3}
\mr{deg}(T^j (\Omega_X) \otimes (\mcF/\mcF^1))&
= \mr{deg} (T^j (\Omega_X)) + \mr{deg}(\mcF/\mcF^1) \cdot \mr{rk}(T^j (\Omega_X))  \\
& = 
\mr{deg}(\Omega_X) \cdot  \binom{\DIM + j-1}{j-1}+ \DEG \cdot  \binom{\DIM + j-1}{j},
 \notag
\end{align}
where $\mr{deg}(-) = c_1 (-) \cdot H^{\DIM -1}$.
Hence, we have
\begin{align} \label{Er443w}
\mr{deg} (\mcF) &= \sum_{j=0}^{n-1} \mr{deg}(\mcF^j/\mcF^{j+1}) \\
&= \sum_{j=0}^{n-1} \mr{deg} (T^j (\Omega_X) \otimes (\mcF/\mcF^1)) \notag \\
&\! \stackrel{(\ref{FG3})}{=} \sum_{j=0}^{n-1} 
\left(\mr{deg}(\Omega_X) \cdot \binom{\DIM + j-1}{j-1}+ \DEG \cdot  \binom{\DIM + j-1}{j}\right) \notag \\
& =\mr{deg}(\Omega_X) \cdot \binom{\DIM+ n-1}{n-2}+ \DEG \cdot \binom{\DIM + n -1}{n-1} \notag \\
& = \left(\mr{deg}(\Omega_X) + \frac{\DEG (\DIM +1)}{n-1} \right) \binom{\DIM + n-1}{n-2}. \notag
\end{align}
On the other hand, according to  ~\cite[Corollary 3.2.4]{LeQu}, 
 the natural morphism $F^{(N)*}_{X/k}(\mcF^\nabla) \migi \mcF$
extending the inclusion $\mcF^\nabla \migiincl \mcF$ is an isomorphism.
This implies that $\mr{deg} (\mcF)$ is equal to $p^N \cdot \mr{deg}(\mcF^\nabla)$ and hence  divisible by  $p^\N$.
By this fact together with (\ref{Er443w}),  the proof of the assertion is completed.
\end{proof}

\LSP
\subsection{}
 \label{SS7ww9}
 Denote by $\mr{Pic}(X^{(N)})$ the group of  line  bundles on $X^{(N)}$.
If
$\mr{Op}_{n, N}^{^\mr{Zzz...}}$ denotes 
the set of isomorphism classes of dormant $(n, N)$-opers 
   on $X$,
  then  we can define  an action of $\mr{Pic}(X^{(N)})$
on $\mr{Op}_{n, N}^{^\mr{Zzz...}}$ as follows:
Let $\mcN$ be a line bundle on $X^{(N)}$ and $\msF^\heartsuit := (\mcF, \nabla, \{ \mcF^j \}_j)$  a dormant $(n, N)$-oper  on $X$.
Denote by $\nabla_{\mcN, \mr{can}}^{(N-1)} \otimes \nabla$ the $\mcD_X^{(N-1)}$-action on the tensor product $F_{X/k}^{(N)*}(\mcN) \otimes \mcF$ induced  by $\nabla_{\mcN, \mr{can}}^{(N-1)}$ and $\nabla$ in a natural manner.
Then, one may verify that the collection of data
\begin{align}
\msF_{\otimes \mcN}^\heartsuit := (F_{X/k}^{(N)*}(\mcN) \otimes \mcF, \nabla_{\mcN, \mr{can}}^{(N-1)} \otimes \nabla, \{ F_{X/k}^{(N)*}(\mcN) \otimes \mcF^j \}_{j=0}^n)
\end{align}
forms a dormant $(n, N)$-oper  on $X$.
The resulting assignment $(\mcN, \msF^\heartsuit) \mapsto \msF^\heartsuit_{\otimes \mcN}$ defines a desired action of $\mr{Pic}(X^{(N)})$ on $\mr{Op}_{n, N}^{^\mr{Zzz...}}$.
In particular, we obtain the quotient set
\begin{align}
\overline{\mr{O}}\mr{p}_{n, N}^{^\mr{Zzz...}} := \mr{Op}_{n, N}^{^\mr{Zzz...}}/\mr{Pic}(X^{(N)}).
\end{align}
We write $[\msF^\heartsuit]$ for the element of  $\overline{\mr{O}}\mr{p}_{n, N}^{^\mr{Zzz...}}$ represented by $\msF^\heartsuit$.
In the  case of $\mr{dim} (X) = 1$,
 this set may be identified with the set of dormant $\mr{PGL}_2$-opers of level $N$ on $X$, in the sense of ~\cite[Definition 4.2.5]{Wak9}.

\SSP
\bpr \label{C100}
Suppose that  $\mr{dim} (X) =1$ (resp.,  $n=2$).
Then, the map of sets 
\begin{align} \label{E500}
\coprod_{d \in \mbZ_{>0}}\mr{Op}_{(n, N, \DEG), +\mr{bir}}^{^\mr{Zzz...}} \migi \overline{\mr{O}}\mr{p}_{n, N}^{^\mr{Zzz...}}
\ \left(\text{resp.,} \ 
\coprod_{d \in \mbZ_{>0}}\mr{Op}_{(2, N, d), +\mr{imm}}^{^\mr{Zzz...}} \migi \overline{\mr{O}}\mr{p}_{2, N}^{^\mr{Zzz...}} \right)
\end{align}
given by assigning   $(\msF^\heartsuit, q) \mapsto [\msF^\heartsuit]$
 is surjective.
\epr
\begin{proof}
We only consider the resp'd assertion since  the non-resp'd assertion can be proved by an entirely similar argument  (by applying assertion (i) of Proposition \ref{L23} instead of  (ii)). 

Let $\msF^\heartsuit := (\mcF, \nabla, \{ \mcF^j \}_j)$ be a dormant $(2, N)$-oper  on  $X$.
Denote by $\pi : \mbP (\mcF^\nabla) \migi X^{(N)}$ the natural projection.
Then, we can find a  very ample  line bundle $\mcN$ on $X^{(N)}$ such that the line bundle $\mcO_{\mbP (\mcF^\nabla)}(1) \otimes \pi^*(\mcN)$ on $\mbP (\mcF^\nabla)$ is very ample.
Write $\sigma_\mbP$ for the closed immersion $\mbP (\mcF^\nabla) \migiincl \mbP^{\DD}$ (where $\DD > 0$) defined by the complete linear system associated to $\pi^*(\mcN) \otimes \mcO_{\mbP (\mcF^\nabla)}(1)$.
Since $\pi_* (\pi^*(\mcN) \otimes \mcO_{\mbP (\mcF^\nabla)}(1))$ may be identified with $\mcN \otimes \mcF^\nabla$ by the projection formula,
$\sigma_\mbP$ coincides with  the morphism induced, via projectivization, 
from the natural  morphism 
\begin{align} \label{E501}
H^0 (X^{(N)}, \mcN \otimes \mcF^\nabla)\otimes_k \mcO_{X^{(N)}}\migi \mcN \otimes \mcF^\nabla.
\end{align}
Hence, the pull-back of (\ref{E501}) by $F_{X/k}^{(N)}$ gives, after  choosing an identification $H^0 (X^{(N)}, \mcN \otimes \mcF^\nabla) = k^{\DD +1}$,   a surjective morphism of $\mcD_X^{(N-1)}$-modules
\begin{align} \label{E502}
q : (\mcO_X, \nabla_{\mr{triv}}^{(N-1)})^{\oplus (\DD +1)} \migisurj (F_{X/k}^{(N)*}(\mcN) \otimes \mcF, \nabla_{\mcN, \mr{can}}^{(N-1)} \otimes \nabla).
\end{align}
Since the composite
$\overline{\sigma} : X \xrightarrow{\sigma} \mbP (\mcF^\nabla) \xrightarrow{\sigma_\mbP} \mbP^{\DD}$ is a closed immersion by Proposition \ref{L23}, (ii), 
the pair $(\msF^\heartsuit_{\otimes \mcN}, q)$ specifies an element of 
$\coprod_{d \in \mbZ_{>0}}\mr{Op}_{(2, N, d), +\mr{imm}}^{^\mr{Zzz...}}$ mapped to $[\msF^\heartsuit] \in \overline{\mr{O}}\mr{p}_{2, N}^{^\mr{Zzz...}}$ via (\ref{E500}).
This implies the surjectivity of (\ref{E500}),  so the proof of this proposition is completed.
\end{proof}

\LSP
\subsection{}
\label{SS71}

Let us  construct, under a certain condition,  a dormant $(n, N)$-oper by using the Gauss map of order $n-1$ on $X$ associated to  a closed immersion
 into a projective space.
To this end, we shall consider the following two conditions on a quadruple  $(X, n, N, \DEG)$:
\begin{itemize}
\item[(a)]
 $X$ is a smooth projective curve over $k$ of genus $g>1$ and $n$, $N$, $\DEG$ are positive integers with $1 < n < p$ and  $p \nmid  (g-1)$;
\item[(b)]
$X$ is a smooth projective variety over $k$ whose  tangent bundle $\mcT_X$
   is stable with respect to some ample line bundle and  $n$, $N$, $\DEG$ are positive integers with $n=2 <p$, $p \nmid d$.
\end{itemize}
Then, the following assertion holds.
\SSP
\bpr \label{L13}
Let $(X, n, N, \DEG)$ be  a quadruple satisfying one of the conditions (a),  (b) described above.
Let  $\mcL$ be a line bundle on $X$ of degree $\DEG$ with respect to some closed immersion from $X$ to a projective space.
Also, let $\nabla$ be 
a left  $\mcD_X^{(\N-1)}$-action on $J_{n-1} (\mcL)$
extending its $\mcO_X$-module structure whose 
 $p$-$(N-1)$-curvature vanishes identically.
Then,  the collection
\begin{align} \label{E81}
(J_{n-1} (\mcL), \nabla, \{ J_{n-1} (\mcL)^j \}_{j=0}^{n})
\end{align}
(cf. (\ref{E70}) for the definition of the filtration $\{ J_{n-1} (\mcL)^j \}_{j}$) forms a dormant $(n, N)$-oper   on $X$.
\epr
\begin{proof}
 First, we shall consider the case where  the condition (a) is satisfied.
Let $B$ be the subset of $\{1, \cdots, n-1 \}$ consisting of integers $j$ satisfying $\nabla^{(0)}(J_{n-1}(\mcL)^j) \subseteq \Omega_X \otimes J_{n-1}(\mcL)^{j-1}$.
Suppose that $B \neq \{1, \cdots, n-1 \}$.
Then,  there exists    the minimum number $j_0$  in $\{ 1, \cdots, n-1 \} \setminus B$.
Since $1 \in B$, we have $j_0 \geq 2$.
The integer  $j_0 -1$ belongs to $ B$,  so  the following $\mcO_X$-linear composite can be defined:
\begin{align} \label{E78}
J_{n-1}(\mcL)^{j_0-1} &\xrightarrow{\nabla^{(0)}} \Omega_X  \otimes J_{n-1}(\mcL)^{j_0-2} \\
&\migisurj \Omega_X  \otimes (J_{n-1}(\mcL)^{j_0-2}/J_{n-1}(\mcL)^{j_0-1}) \left(\stackrel{(\ref{E77})}{\cong} \Omega_X^{\otimes (j_0-1)} \otimes \mcL \right). \notag
\end{align}
It follows from $\mr{deg}(\Omega_X) > 0$ and (\ref{E77}) that
this composite 
becomes the zero map   when restricted to  $J_{n-1}(\mcL)^{j_0} \left(\subseteq  J_{n-1}(\mcL)^{j_0-1}\right)$.
This implies $\nabla^{(0)}(J_{n-1}(\mcL)^{j_0}) \subseteq \Omega_X \otimes J_{n-1}(\mcL)^{j_0-1}$, which contradicts the fact  that $j_0 \notin B$.
Hence,  the equality $B = \{1, \cdots, n-1 \}$ holds. 

Now, let us fix  $j \in \{1, \cdots, n-1\}$, and 
 denote by $\mr{KS}^j$ the $j$-th Kodaira-Spencer map (cf. (\ref{e100})) associated to the collection (\ref{E81}).  
We shall prove the claim that $\mr{KS}^j$ is nonzero.
Suppose, on the contrary, that  $\mr{KS}^{j} =0$.
Then, $J_{n-1} (\mcL)^{j}$ is closed under $\nabla^{(0)}$, and
we can define a connection $\nabla^{(0)}_{j}$ on $J_{n-1} (\mcL)/J_{n-1} (\mcL)^{j}$ induced from  $\nabla^{(0)}$ via the quotient $J_{n-1} (\mcL) \migisurj J_{n-1} (\mcL)/J_{n-1} (\mcL)^{j}$.
Since  $\nabla^{(0)}_{j}$ has vanishing $p$-curvature, 
the $\mcO_X$-linear  morphism 
\begin{align} \label{E438}
F^*_{X/k} (\mr{Ker}(\nabla^{(0)}_{j})) \migi J_{n-1} (\mcL)/J_{n-1} (\mcL)^{j}
\end{align}
 extending the  inclusion $\mr{Ker}(\nabla^{(0)}_{j}) \migiincl J_{n-1} (\mcL)/J_{n-1} (\mcL)^{j}$
(regarded as an $\mcO_{X^{(1)}}$-linear morphism via the underlying homeomorphism of $F_{X/k}$) is an isomorphism (cf. ~\cite[Theorem (5.1)]{Kal}).
By putting  $\nabla_n^{(0)} := \nabla^{(0)}$, we obtain 
the following sequence of equalities for each $j' \in \{ j, n \}$:
\begin{align}
p \cdot \mr{deg}(\mr{Ker}(\nabla_{j'}^{(0)}))&=  \mr{deg} (F_{X/k}^* (\mr{Ker}(\nabla_{j'}^{(0)}))) \\
& = \mr{det} (J_{n-1} (\mcL)/J_{n-1} (\mcL)^{j'})  \notag \\
&= \sum_{i=0}^{j'-1} \mr{deg} (J_{n-1}(\mcL)^i/J_{n-1} (\mcL)^{i+1}) \notag\\ 
& = \sum_{i=0}^{j'-1} \mr{deg} (\Omega_{X}^{\otimes i} \otimes \mcL) \notag \\
& = \sum_{i=0}^{j'-1} \left(i \cdot (2g-2) + \DEG\right) \notag \\
& =  j'  \cdot ((j'-1) \cdot (g-1) + \DEG). \notag
\end{align}
This implies (from the assumption $n < p$) that both $(j-1) \cdot (g-1) + \DEG$ and $(n-1) \cdot (g-1) + \DEG$  are  divisible by $p$.
In particular, the integer
\begin{align}
(g-1)(n-j) \left(=  ((n-1) \cdot (g-1) + \DEG)- ((j-1) \cdot (g-1) + \DEG)\right)
\end{align}
is divisible by $p$.
This contradicts the assumption that $p \nmid (g-1)$.
Hence,  $\mr{KS}^{j}$ turns out to be  nonzero,  and this completes the proof of the claim.

Moreover, by comparing the restrictive degrees of  the line bundles
$J_{n-1}(\mcL)^{j}/J_{n-1}(\mcL)^{j+1}$ and  $\Omega_X  \otimes (J_{n-1}(\mcL)^{j-1})/J_{n-1}(\mcL)^{j})$, 
we see that $\mr{KS}^{j}$ is an isomorphism.
Consequently, the collection (\ref{E81}) forms a dormant $(n, N)$-oper on $X$.

Next,  let us consider the case where the condition (b) is satisfied.
Suppose that the $1$-st Kodaira-Spencer map $\mr{KS}^1 : J_1 (\mcL)^1 \migi \Omega_X \otimes (J_1 (\mcL)/J_1 (\mcL)^1)$ associated to (\ref{E81}) coincides with  the  zero map.
This implies that $J_1 (\mcL)^1$ is closed under $\nabla^{(0)}$, so we can define a connection $\nabla_1^{(0)}$ on $\mcL \left(= J_1 (\mcL)/J_1 (\mcL)^1\right)$ induced naturally from $\nabla^{(0)}$. 
By an argument similar to the above argument,
we have $d \left(= \mr{deg}(\mcL) \right) =p \cdot \mr{deg}(\mr{Ker}(\nabla_1^{(0)}))$, which contradicts the assumption that $p \nmid d$.
Hence, $\mr{KS}^1$ specifies a {\it nonzero} endomorphism of $\Omega_X \otimes \mcL \left(=  J_1 (\mcL)^1 =\Omega_X \otimes (J_1 (\mcL)/J_1 (\mcL)^1)\right)$.
Since $\mcT_X$ (hence also $\Omega_X \otimes \mcL$) is stable,
$\mr{KS}^1$ must be an isomorphism.
That is to say, the collection (\ref{E81}) defines a dormant $(2, N)$-oper.
This completes the proof of the proposition.
 \end{proof}
\SSP

By applying the above proposition, we obtain 
the following assertion, which is the main result of this section.

\SSP
\bt \label{T13}
Let $(X, n, N, \DEG)$ be  a quadruple satisfying one of the conditions (a),  (b).
We shall set $\chi := (n, N, \DEG)$ and  $\lambda:= (n-1, N, d)$.
Then, the following assertions hold:
 \begin{itemize}
\item[(i)]
Let $\iota : X \migiincl  \mbP^{\DD}$ (where  $\DD >0$)  be a closed immersion classified by 
$\mr{G au}_{\lambda}^\mr{F}$.
Write $\mcL := \iota^* (\mcO_{\mbP^\DD}(1))$, and 
we shall  set
\begin{align} \label{e9899}
\msF_\iota^\heartsuit := (J_{n-1}(\mcL), \nabla_{\iota, \mr{Gau}}^{(N-1)}, \{ J_{n-1}(\mcL)^j \}_{j=0}^n).
\end{align}
Then, the pair $\mpf_\iota := (\msF_\iota^{\heartsuit}, \alpha_\iota^{n-1})$ (cf. (\ref{E50}) for the definition of $\alpha_\iota^{n-1}$) 
specifies an element of  $\mr{Op}_{\chi, +\mr{imm}}^{^\mr{Zzz...}}$.
Moreover,  the map of sets 
\begin{align} \label{E83}
\Xi_{\chi} : \mr{G au}_{\lambda}^\mr{F} \migi
\mr{Op}_{\chi, +\mr{imm}}^{^\mr{Zzz...}}.
\end{align}
given by  $\iota \mapsto \mpf_\iota$ is injective.
\item[(ii)]
Let $N'$ be a positive integer with $N' < N$.
We shall set $\chi' := (n, N', \DEG)$ and $\lambda' := (n-1, N', \DEG)$.
Then, the following square diagram is commutative:
\begin{align} \label{E02dd}
\vcenter{\xymatrix@C=46pt@R=36pt{
\mr{G au}_{\lambda}^\mr{F} \ar[r]^-{\Xi_{\chi}} \ar[d]_-{\mr{inclusion}} & \mr{Op}_{\chi, +\mr{imm}}^{^\mr{Zzz...}} \ar[d]^-{\mr{inclusion}}
\\
\mr{G au}_{\lambda'}^\mr{F} \ar[r]_-{\Xi_{\chi'}} & \mr{Op}_{\chi', +\mr{imm}}^{^\mr{Zzz...}}.
}}
\end{align}

\item[(iii)]
Suppose further that $n=2$.
Then, the map $\Xi_{\chi}$ is bijective.
\end{itemize}
\et
\begin{proof}
First, we shall consider assertion (i).
It follows from Proposition \ref{L13} that $\msF_\iota^\heartsuit$ forms a dormant $(n, N)$-oper  on $X$.
Since $\gamma_\iota^{n-1}$ factors through $F_{X/k}^{(N)}$,
there exists a morphism $h : X^{(N)} \migi \mr{Grass}(\binom{\DIM + n-1}{n-1}, \DD +1)$ with $h \circ F_{X/k}^{(N)} = \gamma_\iota^{n-1}$.
If  $q_0 : \mcO_{X^{(N)}}^{\oplus (\DD +1)} \migisurj \mcQ$ denotes the pull-back of $q_{\mr{univ}}$ (cf. (\ref{e37})) by $h$, then
$F_{X/k}^{(N)*}(q_0)$ may be identified with $\alpha_\iota^{n-1}$.
This implies that $\alpha_\iota^{n-1}$ defines a surjection   of $\mcD_X^{(N-1)}$-modules $(\mcO_X, \nabla_{\mr{triv}}^{(N-1)})^{\oplus (\DD +1)} \migisurj (J_{n-1}(\mcL), \nabla_{\iota, \mr{Gau}}^{(N-1)})$.
Thus, the pair $(\msF_\iota^\heartsuit, \alpha_\iota^{n-1})$   turns out  to be an element of $\mr{Op}_{\chi, +\mr{imm}}^{^\mr{Zzz...}}$.
Moreover, 
the injectivity  of $\Xi_{\chi}$ follows immediately  from the observation that
each  closed immersion $\iota : X \migiincl \mbP^{\DD}$ in $\mr{G au}_{\lambda}^\mr{F}$ may be reconstructed as the projectvization of  the composite of $\alpha_\iota^{n-1}$ and the natural quotient $J_{n-1}(\iota^*(\mcO_{\mbP^{\DD}}(1))) \migisurj \iota^*(\mcO_{\mbP^{\DD}}(1))$.
This completes the proof of assertion (i).

Also, assertion (ii) follows from the definition of $\Xi_\chi$.

Next, we shall prove assertion (iii), i.e., the  surjectivity of $\Xi_{\chi}$ under the assumption that $n=2$.
Let $\mpf := (\msF^\heartsuit, q)$ (where $\msF^\heartsuit := (\mcF, \nabla, \{ \mcF^j \}_{j=0}^2)$)
be a pair classified by $\mr{Op}_{\chi, +\mr{imm}}^{^\mr{Zzz...}}$.
It follows from Proposition  \ref{L23}, (ii), that the morphism $\iota := \iota_{\mpf} : X \migi \mbP^{\DD_\mpf}$ (cf. (\ref{E4444})) is a closed immersion.
Let us write $\mcL := \iota^* (\mcO_{\mbP^{\DD_\mpf}}(1))$.
Also, write $\sigma : X \migi \mbP (\mcF)$ (resp., $\overline{\sigma} : X \migi \mbP (\mcF^\nabla)$) for the morphism induced from the surjection $\mcF \migisurj \mcF/\mcF^1$ (resp., the composite $\mcF^\nabla \migiincl \mcF \migisurj \mcF/\mcF^1$) as defined  in Proposition \ref{L23}.  
Then, the  surjection
$\iota^* (\Omega_{\mbP^{\DD}}) \migisurj \Omega_{X}$ induced by $\iota$ can be decomposed  as
$\iota^*(\Omega_{\mbP^{\DD}}) \migisurj  \overline{\sigma}^*(\Omega_{\mbP (\mcF^\nabla)}) \migisurj \Omega_{X}$.
Since the differential  of the composite $X \xrightarrow{\overline{\sigma}} \mbP(\mcF^\nabla) \xrightarrow{\mr{projection}}X^{(N)}$ (which coincides with $F_{X/k}^{(N)}$) is the zero map,
the surjection  $\overline{\sigma}^*(\Omega_{\mbP (\mcF^\nabla)}) \migisurj \Omega_{X}$  factors through
the   quotient
$\overline{\sigma}^*(\Omega_{\mbP (\mcF^\nabla)})\migisurj \overline{\sigma}^*(\Omega_{\mbP (\mcF^\nabla)/X^{(N)}})$.
The resulting morphism  between line bundles $\left(\sigma^*(\Omega_{\mbP (\mcF)/X}) =  \right) \overline{\sigma}^*(\Omega_{\mbP (\mcF^\nabla)/X^{(N)}}) \migi \Omega_X$ is surjective, hence it is also  bijective.
This implies that  
the families  of linear subvarieties  in $\mbP^{\DD_\mpf}$ (parametrized by $X$)  given by $q$ and $\alpha_\iota^{n-1}$, respectively,  are identical, i.e., $\mbP (J_{1}(\mcL)) = \mbP (\mcF)$.
 It follows from 
 the various definitions involved that
the associated  isomorphism  $J_{n-1}(\mcL) \isom \mcF$ between quotient bundles of  $\mcO_X^{\oplus (\DD_\mpf +1)}$
defines  an isomorphism 
  $\mpf \cong (\msF_\iota^\heartsuit, \alpha_\iota^{n-1})$.
  This shows  the surjectivity of $\Xi_{\chi}$, so
 the proof of assertion (iii) is  completed.
\end{proof}

\vspace{10mm}
\section{Purely inseparable Gauss maps on a curve} \label{S08}
\SSP

In this section,  we 
consider  a sufficient condition for the nonemptiness of the set $\mr{Op}_{\chi, + \Box}^{^\mr{Zzz...}}$ (where $\Box \in \{\mr{bir}, \mr{imm} \}$) in the case of  $\mr{dim} (X) =1$.
As an application of this result, 
we shows (cf. Theorem \ref{C091}) that, for any $N>0$,  there always exists a closed immersion $X \migiincl \mbP^\DD$ with purely inseparable Gauss map of degree 
    $p^N$.

\LSP
\subsection{}
 \label{SS11}

Let $\chi := (n, N, \DEG)$  be a triple of positive  integers with $n>1$, and let
 $X$ be  a smooth projective curve over $k$ of genus $g>1$.

\SSP
\bpr \label{P8}
Let $\DD$ be  a positive  integer.
 Suppose that 
   there exists an integer $a$ satisfying
 \begin{align} \label{E54}
 \frac{\DD +1}{n} + g-1 \geq \frac{d +(g-1)(n-1)}{p^N} = a \geq \frac{(g-1)(n-1)}{p^N} + 2g+1.
 \end{align}
Also, suppose that $n<p$ (resp., $n=2 <p$).
 Then, 
there exists an element  $\mpf := (\msF^\heartsuit, q)$ of  $\mr{Op}_{\chi, + \mr{bir}}^{^\mr{Zzz...}}$ (resp., $\mr{Op}_{\chi, + \mr{imm}}^{^\mr{Zzz...}}$) 
such that $\DD_\mpf = \DD$ and the underlying  vector bundle of $\msF^\heartsuit$ has degree $p^N \cdot n \cdot a$.
In particular,  if $\DEG$  is   sufficiently large  relative to $g$, $N$ and  satisfies  $d \equiv -(g-1)(n-1)$ mod $p^N$, 
then the set  $\mr{Op}_{\chi, + \mr{bir}}^{^\mr{Zzz...}}$ (resp., $\mr{Op}_{\chi, + \mr{imm}}^{^\mr{Zzz...}}$) is nonempty.
\epr
\begin{proof}
We only consider the non-resp'd assertion since the resp'd assertion can be proved by an entirely similar argument (by applying assertion (ii) of Proposition \ref{L23} instead of  (i).)

Let us take a theta characteristic  of $X$, i.e., a line bundle $\varTheta$ on $X$ together with an isomorphism $\varTheta^{\otimes 2} \isom \Omega_X$.
According to ~\cite[Theorem 7.5.2]{Wak6},
there exists a dormant $(2, N)$-oper $\msF^\heartsuit := (\mcF_0, \nabla_0, \{ \mcF_0^j \}_{j=0}^2)$  on $X$ with 
$\mcF_0^1 = \varTheta$ and $\mcF^0_0/\mcF_0^1 = \varTheta^\vee$.
Denote by $S^{n-1}(\mcF_0)$ the $(n-1)$-st symmetric product of $\mcF_0$ over $\mcO_X$.
Note that  $S^{n-1}(\mcF_0)$ forms a rank $n$ vector bundle on $X$ of degree $0$  and admits a $\mcD_X^{(N-1)}$-action $S^{n-1}(\nabla_0)$  induced  naturally by $\nabla_0$.
Moreover,  $S^{n-1}(\mcF_0)$ is equipped  with an $n$-step decreasing filtration $\{ S^{n-1}(\mcF_0)^j\}_{j=0}^n$ induced from $\{ \mcF_0^j \}_j$;
to be precise,  we set $S^{n-1}(\mcF_0)^0 := S^{n-1}(\mcF_0)$, $S^{n-1}(\mcF_0)^n := 0$, and $S^{n-1}(\mcF_0)^j$ (for each $j =1, \cdots , n-1$) is defined as  
the image of $(\mcF_0^1)^{\otimes j} \otimes \mcF_0^{\otimes (n-j-1)}$ via the natural quotient $\mcF_0^{\otimes (n-1)} \migisurj S^{n-1}(\mcF_0)$.
This filtration satisfies that 
\begin{align} \label{E61}
S^{n-1}(\mcF_0)^j/S^{n-1}(\mcF_0)^{j+1} \cong \varTheta^{\otimes (1-n)} \otimes \Omega_X^{\otimes j}
\end{align}
 for every $j=0, \cdots, n-1$.
 Since $\mr{KS}_{\msF_0^\heartsuit}^1$ is an isomorphism, 
  the assumption  $n < p$ implies that the Kodaira-Spencer maps associated to  the collection 
\begin{align}
(S^{n-1}(\mcF_0), S^{n-1}(\nabla_0), \{ S^{n-1}(\mcF_0)^j \}_{j=0}^n)
\end{align}
are verified to be isomorphisms.
That is to say, this collection
 forms a dormant $(n, N)$-oper  on $X$.

Since $\nabla_0$ (hence also $S^{n-1}(\nabla_0)$) has vanishing $p$-$(N-1)$-curvature,
the inclusion  $S^{n-1}(\mcF_0)^\nabla \migiincl S^{n-1}(\mcF_0)$ extends to 
an isomorphism of $\mcD_X^{(N-1)}$-modules  
\begin{align} \label{E43}
(F^{(N)*}_{X/k}(S^{n-1}(\mcF_0)^\nabla), \nabla^{(N-1)}_{S^{n-1}(\mcF_0)^\nabla, \mr{can}}) \isom (S^{n-1}(\mcF_0), S^{n-1}(\nabla_0))
\end{align}
(cf. ~\cite[Corollary 3.2.4]{LeQu}).
Hence, the faithful flatness of $F_{X/k}^{(N)}$ implies that  $S^{n-1}(\mcF_0)^\nabla$ forms a rank $n$ vector bundle of degree $0 \left(= \frac{1}{p^N} \cdot \mr{deg}(S^{n-1}(\mcF_0)) \right)$ on $X^{(N)}$.

Now, let us choose a line bundle $\mcN$ on $X^{(N)}$  of degree $a$ and
 a quotient line bundle $\mcM$ of $\mcN \otimes S^{n-1}(\mcF_0)^\nabla$ having minimal degree.
Write 
\begin{align}
\mcF := F^{(N)*}_{X/k}(\mcN) \otimes S^{n-1}(\mcF_0) \ \  \text{and} \ \  \mcF^j := F^{(N)*}_{X/k}(\mcN) \otimes S^{n-1}(\mcF_0)^j
\end{align}
 ($j=0, \cdots, n$).
 The degree of $\mcF$ is given by
 \begin{align} \label{ER443}
 \mr{deg}(\mcF) = n \cdot \mr{deg}(F^{(N)*}_{X/k}(\mcN)) + \mr{deg}( S^{n-1}(\mcF_0)) = n \cdot p^N \cdot \mr{deg}(\mcN) + 0 = p^N \cdot n \cdot a.
 \end{align}
 Also, we have 
 \begin{align} \label{E4446}
 \mr{deg}(\mcF/\mcF^1) = \mr{deg}(F_{X/k}^{(N)*}(\mcN) \otimes \varTheta^{\otimes (1-n)}) = p^N  \cdot a + (g-1)(1-n) = d.
 \end{align}
By (\ref{E61}) and (\ref{E43}),
$F_{X/k}^{(N)*}((\mcN \otimes S^{n-1}(\mcF_0)^\nabla)^\vee)$ may be identified with
$\mcF^\vee$
   and has a filtration whose graded pieces  are  isomorphic to the line bundles
 \begin{align} \label{E70}
   F^{(N)*}_{X/k}(\mcN^\vee) \otimes \varTheta^{\otimes (n-1)} \otimes \Omega_{X}^{\otimes (-j)}
 \end{align}
 ($j=0, \cdots, n-1$).
 Hence, since $F^{(N)*}_{X/k}(\mcM^\vee)$ specifies  a line subbundle of  $F_{X/k}^{(N)*}((\mcN \otimes S^{n-1}(\mcF_0)^\nabla)^\vee)$,
  we have
 \begin{align} \label{E4445}
 \mr{deg}(\mcM) & = -\frac{1}{p^N} \cdot \mr{deg}(F^{(N)*}_{X/k}(\mcM^\vee)) \\
 &\geq -\frac{1}{p^N} \cdot \mr{max}\left\{ \mr{deg}(F^{(N)*}_{X/k}(\mcN^\vee) \otimes \varTheta^{\otimes (n-1)} \otimes \Omega_X^{\otimes (-j)}) \, \Big| \, 0 \leq j \leq n-1\right\}\notag \\
 &= -\frac{1}{p^N} \cdot \mr{deg}(F^{(N)*}_{X/k}(\mcN^\vee) \otimes \varTheta^{\otimes (n-1)}) \notag \\
 &=  a - \frac{1}{p^N} \cdot  (n-1) (g-1) \notag \\
 & \geq 2g+1,\notag
 \end{align}
 where the last inequality follows from the assumption (\ref{E54}).
 This implies that $\mcN \otimes S^{n-1}(\mcF_0)^\nabla$ is globally generated and  very ample (cf. ~\cite[Proposition 2, (iii) and (iv)]{IT}), and moreover,  the following equalities hold:
  \begin{align}
 h^1 (\mcN \otimes S^{n-1}(\mcF_0)^\nabla) = h^0 ((\mcN \otimes S^{n-1}(\mcF_0)^\nabla)^\vee \otimes \Omega_X)=0.
 \end{align}
 By the Riemann-Roch theorem, we have 
 \begin{align}
 h^0 (\mcN \otimes S^{n-1}(\mcF_0)^\nabla) &= h^0 (\mcN \otimes S^{n-1}(\mcF_0)^\nabla) - h^1 (\mcN \otimes S^{n-1}(\mcF_0)^\nabla)  \\
 & = n a + n (1-g) \notag \\
 &\leq \DD +1, \notag
 \end{align}
 where the last inequality follows from the assumption (\ref{E54}).
 Hence, there exists an $\mcO_{X^{(N)}}$-linear surjection 
 \begin{align}
 q_0 : \mcO_{X^{(N)}}^{\oplus (\DD +1)} \migisurj \mcN \otimes S^{n-1}(\mcF_0)^\nabla
 \end{align}
  such that the  associated morphism 
 \begin{align} \label{e4601}
\mbP (\mcN \otimes S^{n-1}(\mcF_0)^\nabla) \migi \mbP^{\DD}
\end{align}
   is a closed immersion.
It follows from Proposition  \ref{L23}, (i),  that the following composite is birational onto its image:
\begin{align}
\iota : X 
\migi
 \mbP (\mcF) \left(= X \times_{X^{(N)}} \mbP (\mcN \otimes S^{n-1}(\mcF_0)^\nabla)\right) 
 \xrightarrow{\mr{projection}}\mbP (\mcN \otimes S^{n-1}(\mcF_0)^\nabla)  
 \xrightarrow{(\ref{e4601})} \mbP^{\DD},
\end{align}
 where the first arrow denotes the morphism arising  from 
 the quotient $\mcF \migisurj \mcF/\mcF^1$.
If  $\nabla$  denotes the $\mcD_X^{(N-1)}$-action on $\mcF$ induced by $\nabla_{\mcN, \mr{can}}^{(N-1)}$ and $S^{n-1}(\nabla_0)$, then
  the  collection of data
\begin{align} \label{e14}
\msF^\heartsuit := (\mcF, \nabla, \{ \mcF^j \}_{j=0}^n),
\end{align}
 forms a dormant $(n, N)$-oper on $X$ (cf. ~\cite[\S\,4.2]{Wak9}).
Moreover,  the pull-back of $q_0$ by $F_{X/k}^{(N)}$ defines, via (\ref{E43}),  
a surjective morphism of $\mcD_X^{(N-1)}$-modules 
\begin{align}
q: (\mcO_X, \nabla^{(N-1)}_{\mr{triv}})^{\oplus (\DD +1)} \migisurj (\mcF, \nabla).
\end{align}
It follows that the   pair $(\msF^\heartsuit, q)$ 
specifies  an element of $\mr{Op}_{\chi, +\mr{bir}}^{^\mr{Zzz...}}$.
This completes the proof of the assertion.
\end{proof}
\SSP

\bco \label{C09}
Suppose that $2  < p$ and  $p \nmid (g-1)$.
Also, let $\DEG$ be an integer satisfying $p^N (2g+1) \leq d$ 
and $p^N \mid (d + g-1)$.
Then, $\mr{G au}_{(1, N, \DEG)}^\mr{F}$ is nonempty.
If, moreover, the integer $\DEG$ satisfies $p^{N'} \nmid (d + g-1)$ for a positive integer $N' \left(>N\right)$,
then $\mr{G au}_{(1, N', \DEG)}^\mr{F}$ is empty.
\eco
\begin{proof}
The  assertion follows immediately from the resp'd portion of Proposition \ref{P8}  and Theorem \ref{T13}, (iii), in the case where $(X, n, N, \DEG)$ satisfies the condition (a).
\end{proof}

\LSP
\subsection{}
 \label{SggS03}
 Let us describe an assertion improving a result by H. Kaji (cf. Introduction).
We shall denote by 
$K(X)$ the function field of $X$ and by 
$\mcK$ the set of subfields $K$ of $K(X)$ satisfying the following condition:  There exists a closed immersion $\iota$ from $X$ to some projective space such that 
the extension of function fields $K(X)/K (\mr{Im}(\gamma_\iota^1))$ defined by the $1$-st order (i.e., classical)  Gauss map $\gamma_\iota^1$ associated to $(X, \iota)$ coincides with $K(X)/K$.

\SSP
\bt[= Theorem \ref{TB}] \label{C091}
Let $X$ be a smooth projective curve over $k$ of genus $g>1$.
Suppose that $2 <p$, $p\nmid (g-1)$.
Then, the  following equality of sets holds:
\begin{align}
\mcK = \left\{ K(X)^{p^N} \, \Big| \, N \geq 0 \right\},
\end{align}
where $K(X)^{p^N} := \{ v^{p^N} \, | \, v \in K (X) \}$.
\et
\begin{proof}
By ~\cite[Corollaries 2.3 and 4.4]{Kaj2} (cf. the discussion  preceding Theorem \ref{TB}) together with  
 the fact mentioned in
Remark \ref{Rf99},
the problem is reduced to proving  that, for every positive integer $N$, there exists $\DEG >0$ with  $\mr{G au}^\mr{F}_{1, N, d} \setminus \mr{G au}^\mr{F}_{1, N+1, d} \neq \emptyset$.
(In fact, the extension of function fields associated to a closed immersion in $\mr{G au}^\mr{F}_{1, N, d} \setminus \mr{G au}^\mr{F}_{1, N+1, d}$ must be equal to  $K(X)/K(X)^{p^N}$.)
However, 
we can always find an integer $\DEG$ with $p^N (2g+1) \leq \DEG$, $p^N \mid (\DEG + g-1)$, and $p^{N+1} \nmid (\DEG + g-1)$, and Corollary \ref{C09}  implies that such an integer $\DEG$ satisfies  the required  condition.
\end{proof}

\vspace{10mm}
\section{A Frobenius-projective structure on a Fermat hypersurface} \label{S08}
\SSP

In this final section, 
we construct a Frobenius-projective structure on a Fermat hypersurface by applying Theorem \ref{T13} and the previous study of Gauss maps in positive characteristic.
We also show  that this Frobenius-projective structure cannot be lifted to sufficiently large levels.

\LSP
\subsection{}
 \label{SS03}

Let $N$ be a positive integer and $X$   a smooth projective variety over $k$ of dimension $\DIM >0$.
Denote by 
 $\mr{PGL}_{\DIM+1}$
 the projective linear group
  over $k$ of rank $\DIM+1$, which 
  can be identified with the automorphism group  of $\mbP^\DIM$.
Let us denote by $(\mr{PGL}_{\DIM+1})_{\Y}^{(N)}$ the Zariski sheaf of groups on $\Y$
given by $U \mapsto \mr{PGL}_{\DIM+1}(U^{(N)})$ for each open subscheme $U$ of $X$.
Also, denote by $\mcP^{\text{\'{e}t}}_X$  
 the Zariski  sheaf of sets on $X$ that assigns, to each open subscheme  $U$ of $X$, the set of \'{e}tale morphisms $U \migi \mbP^\DIM$.
Note that the sheaf $\mcP^{\text{\'{e}t}}_X$ has a natural   $(\mr{PGL}_{\DIM+1})_X^{(N)}$-action (cf. ~\cite[\S\,1.2]{Wak6}).

Recall that a subsheaf $\mcS^\blacklozenge$ of $\mcP^{\text{\'{e}t}}_X$ is said to be 
a {\bf Frobenius-projective structure of level $N$} (or, {\bf $F^N$-projective structure}, for short) on $X$ if
it is closed under the  $(\mr{PGL}_{\DIM+1})_X^{(N)}$-action on $\mcP^{\text{\'{e}t}}_X$ 
 and forms a $(\mr{PGL}_{\DIM+1})_X^{(N)}$-torsor 
 with respect to the resulting $(\mr{PGL}_{\DIM+1})_X^{(N)}$-action 
   on $\mcS^\blacklozenge$ (cf. ~\cite[Definition 2.1]{Hos2}, ~\cite[Definition 1.2.1]{Wak6}).

According to  ~\cite[Theorem A]{Wak6},
there exists an assignment from a 
dormant $(2, N)$-oper  on $X$ (i.e., a dormant indigenous $\mcD_X^{(N-1)}$-modules, in the sense of ~\cite[Definitions 2.3.1 and 3.2.1]{Wak6})
 to an $F^N$-projective structure (cf. Remark \ref{R90}); moreover, if $p \nmid (\DIM+1)$, then
this assignment defines 
 a bijective correspondence between the set of
  $F^N$-projective structures on $X$ and the set of 
  equivalence classes of dormant $(2, N)$-opers. 
  (We here omit the details of the {\it equivalence relation} on dormant $(2, N)$-opers. When $\mr{dim}(X)=1$, each such equivalence class was  referred, in ~\cite[Definition 4.2.5]{Wak9}, to as a {\it dormant $\mr{PGL}_2$-opers of level $N$}.)

\LSP
\subsection{}
\label{SS035}

As an application of  Theorem \ref{T13},
we can construct an $F^N$-projective structure by using the Gauss map of  a certain  Fermat hypersurface (cf. ~\cite[\S\,7]{Wal}).

Hereinafter, let us fix an integer $\DD >1$, and  suppose that $X$ is  
 the Fermat hypersurface of degree $p^N+1$ in the projective space $\mbP^{\DD}$,
  i.e., the smooth hypersurface
defined by 
 the homogenous polynomial 
\begin{align} \label{E481}
f_N := 
t_0^{p^N+1} + \cdots + t_\DD^{p^N+1}.
\end{align}
Write $\iota : X \migiincl \mbP^\DD$ for  the natural closed immersion.
Let us identify $\mr{Grass} (\DD, \DD+1)$ with $\mbP^\DD$ in such a way that 
 if an $\DD$-plane in $\mbP^\DD$ (i.e., a point of $\mr{Grass} (\DD, \DD+1)$) is  given  by  an  equation $\sum_{i=0}^\DD v_i  \cdot t_i = 0$ ($v_0, \cdots, v_\DD \in k$, $(v_0, \cdots, v_\DD) \neq (0, \cdots, 0)$), then it corresponds to the point $[v_0: \cdots : v_\DD]$ of $\mbP^\DD$.
Under this identification, 
the $1$-st order Gauss map $\gamma_{\iota}^{1} : X \migi \mbP^\DD$  associated to $(X, \iota)$
can be described as the assignment 
\begin{align}
a := [a_0: \cdots : a_\DD] \mapsto \left(\left[\frac{\partial f_N}{\partial t_0} (a): \cdots :\frac{\partial f_N}{\partial t_\DD}(a)\right] = \right) [a_0^{p^N}: \cdots : a_\DD^{p^N}].
\end{align}
That is to say, the morphism $X \migi \mr{Im}(\gamma_\iota^1)$ induced by $\gamma_\iota^1$
  coincides with the $N$-th relative Frobenius morphism $F_{X/k}^{(N)}$ of $X$.
It follows that the closed immersion $\iota$ defines an element of $\mr{G au}_{(1, N,  p^{N}+1)}^{\mr{F}}$.

\SSP
\bt \label{C01}
Suppose that $\DD \neq 3$ and $p >2$.
Then,  
the Fermat hypersurface $X$ of degree $p^N +1$ in $\mbP^\DD$
admits
  a canonical $F^N$-projective structure
  \begin{align} \label{e9999}
  \mcS^\blacklozenge_{\mr{Gau}},
  \end{align}
 which  corresponds to the dormant $(2, N)$-oper $\msF_\iota^\heartsuit$ (cf. (\ref{e9899})).
Moreover, if $p \nmid L (L+1)$,
then $X$ admits no $F^{2N+1}$-projective structures.
\et
\begin{proof}
By the assumption $\DD \neq 3$, 
the tangent bundle $\mcT_X$ is  stable (cf. ~\cite[Remark 3.2]{MT} or ~\cite[Corollary 0.3]{PW}).
In particular,  the quadruple $(X, 2, N, p^N +1)$ satisfies 
 the condition (b) described in \S\,\ref{SS71}.
 It follows that we can define the map $\Xi_{(2, N, p^N+1)}$ asserted  in Theorem \ref{T13}, (i), 
 and the image via this map of the element  $\iota \in \mr{G au}_{(1, N, p^{N}+1)}^{\mr{F}}$ determines a dormant $(2, N)$-oper, or equivalently, an $F^N$-projective structure, on $X$.
This completes the proof of  the former assertion.

Next, let us consider the latter assertion.
For each vector bundle $\mcV$ and an integer $m>0$, we shall use the notation ``$c_m^{\mr{crys}}(\mcV)$"  to denote  the $m$-th crystalline Chern class of $\mcV$, which is an element of the $2m$-th crystalline cohomology group $H_{\mr{crys}}^{2m}(X/W)$ ($W$ denotes the ring of Witt vectors over $k$). 
Denote by $H$  the restriction to $X$  of $c_1^\mr{crys} (\mcO_{\mbP^\DD} (1))$.
Then, the Chern polynomial  $c_t^\mr{crys} (\mcT_X)$ of $\mcT_X$ is  given by
\begin{align} \label{E3321}
c_t^\mr{crys} (\mcT_X) &= (1 + H t)^{\DD+1} (1 + (p^N+1)H t)^{-1}\\
& = (1 + H t)^{\DD+1} (1 - (p^N+1)Ht + ((p^N+1) H t)^2 - ((p^N+1) H t)^3+ \cdots), \notag
\end{align}
where the first equality follows from 
 the natural short exact sequence
\begin{align}
0 \longmigi \mcT_X \longmigi \iota^*(\mcT_{\mbP^\DD}) \longmigi \mcO_X (p^{N}+1) \longmigi 0
\end{align}
and the Euler sequence on $\mbP^\DD$.
Since  $H^4_{\mr{crys}}(X/W) = \{ aH^2 \, | \, a \in W \}$ (cf. ~\cite[Exp.\,XI, Theorem 1.5]{Del}, ~\cite[Chap.\,VII, Remark 1.1.11]{PBer0}, ~\cite[Chap.\,II, Corollary 3.5]{Ill}), the equalities   (\ref{E3321}) implies 
\begin{align}
&  \ \ \ \, c_2^\mr{crys} (\mcT_X) - \frac{1}{\DD^2} \cdot \binom{\DD}{2} \cdot c_1^\mr{crys} (\mcT_X)^2 \\
& = \left(p^{2N} - p^N\DD  + p^N + \frac{\DD^2-\DD}{2} \right) H^2 - \left((\DD-p)^2 \cdot \frac{\DD -1}{2 \DD} \right) H^2 \notag \\
& = \frac{p^{2N} \cdot (L+1)}{2L}  H^2 \notag \\
& \not\equiv 0 \ (\text{mod} \ p^{2N+1}), \notag
\end{align}
where the last ``$\not\equiv$" follows from 
  the assumption $p \nmid L (L+1)$.
  Thus, the assertion follows from ~\cite[Theorem 3.7.1]{Wak6}.
\end{proof}
\SSP

\begin{rema} \label{R059}
In the case of  $\DD =2$, the  variety  $X$ defined by (\ref{E481}) is  known as a {\it Hermitian curve}; this is a smooth projective curve   of  genus $\frac{p^N (p^{N}-1)}{2}$ having  large automorphism groups so that it violates  the classical Hurwitz bound (i.e., $\sharp (\mr{Aut}(X)) \leq 84 (g-1)$).
In fact,  let us define $\mr{U}_3 (p^{2N})$ as  the subgroup of $\mr{GL}_3 (\mbF_{p^{2N}}) \left(\subseteq \mr{GL}_3 (k) \right)$ leaving (\ref{E481}) invariant, and $\mr{PGU}_3 (p^{2N})$ as the
factor of $U_3 (p^{2N})$ modulo its center.  
Then, $\mr{PGU}_3 (p^{2N})$ coincides with the full automorphism group of $X$, and its order is given by  $p^{3N}(p^{3N}+1)(p^{2N}-1) \left(> 84 (\frac{p^N (p^N -1)}{2}-1) \right)$.

Since the closed immersion  $\iota : X \migiincl  \mbP^2$ is compatible with the respective $\mr{PGU}_3 (p^{2N})$-actions on $X$ and $\mbP^2$,
the Gauss $\gamma_\iota^1 : X \migi \mr{Grass}(2, 3)$ is compatible with the respective $\mr{PGU}_3 (p^{2N})$-actions.
By the definition of $\Xi_{(2, N, p^N+1)}$,  the dormant $(2, N)$-oper  obtained from $\iota$, hence also the $F^N$-projective structure $\mcS_{\mr{Gau}}^\blacklozenge$,   turns out to be  invariant under the $\mr{PGU}_3 (p^{2N})$-action on $X$.
 So $\mcS_{\mr{Gau}}^\blacklozenge$  has large symmetry in this sense and descends to any \'{e}tale quotient of $X$.
\end{rema}

\SSP

\begin{rema} \label{R05}
According to ~\cite[Corollary]{Shim} (or, ~\cite[Proposition 1]{Shio}, ~\cite[Theorem III]{ShiKa}),
the Fermat hypersurface $X \left(\subseteq \mbP^\DD \right)$ is unirational when $\DD >3$.  
Recall that if any unirational projective variety over the field of complex numbers $\mbC$ 
admits no  projective structures unless it is isomorphic to a projective space; this is  because such a variety 
 contains a rational curve (cf. ~\cite[Theorem 4.1]{JR1}).
In this sense,
the example of  an $F^N$-projective structure resulting from the above theorem
embodies  an exotic phenomenon of algebraic geometry in positive  characteristic.
\end{rema}

\SSP

\begin{rema} \label{R04}
Recall that there is an ``affine" version of an $F^N$-projective structure, which is called an {\it $F^N$-affine structure} (cf. ~\cite[Definition 2.1]{Hos3}, ~\cite[Definition 1.2.1]{Wak6}).
By  a change of structure group from the group  of affine transformations to  
that of projective transformations,
 each  $F^N$-affine structure yields an $F^N$-projective structure.
The only previous examples of $F^N$-projective structures on  higher-dimensional  varieties except for those on projective spaces   were obtained, in that manner,  from  $F^N$-affine projective structures on Abelian varieties or  smooth curves equipped with  a Tango structure via, e.g.,  taking products, \'{e}tale coverings, or quotients by a finite group action (cf. ~\cite[\S\,6.1, \S\,6.5, \S\,8.1]{Wak6}).
On the other hand,
the degree $p^N +1$ Fermat hypersurface $X$ embedded in $\mbP^\DD$   with  $p \nmid \DD$ satisfies
 $c_1^{\mr{crys}} (\mcT_X)\left(  = (\DD -p^N)H\right) \not\equiv 0$ (mod $p^N$) in $H^2_\mr{crys}(X/W) \left(= \{ a H \, | \, a\in W \}\right)$.
Hence, it follows from    ~\cite[Theorem 3.7.1]{Wak6} that $X$ admits no $F^N$-affine structures.
In particular,
 the $F^N$-projective structure $\mcS^\blacklozenge_{\mr{Gau}}$ resulting from the above theorem does not come from  any $F^N$-affine structure via changing the structure group.
This means that 
$\mcS^\blacklozenge_{\mr{Gau}}$ is essentially a new example constructed in a way that has never been done before.
\end{rema}

\end{document}